\documentclass[a4paper,12pt]{article}

\usepackage{amsmath}
\usepackage{amssymb}
\usepackage{amsthm}
\usepackage{mathtools}
\usepackage{enumitem}
\usepackage{empheq}
\usepackage{authblk}

\usepackage[final]{hyperref}
\hypersetup{
	colorlinks=true,
	citecolor=black,
	filecolor=black,
	linkcolor=blue,
	urlcolor=blue
}

\newtheorem{theorem}{Theorem}[section]
\newtheorem{lemma}[theorem]{Lemma}
\newtheorem{proposition}[theorem]{Proposition}

\newtheorem{corollary}[theorem]{Corollary}

\theoremstyle{remark}

\def\he{\hat{e}}
\def\hr{\hat{r}}
\def\hth{\hat{\theta}}
\def\hvp{\hat{\varphi}}
\def\R{\mathbb{R}}
\def\d{\partial}
\def\dd{\textnormal{\textrm{div}}}
\def\dif{{\mathrm d}}
\def\tr{\tilde{\rho}}
\def\tu{\tilde{u}}
\def\tU{\tilde{U}}
\def\tf{\tilde{f}}
\def\br{\rho_{+}}
\def\bq{q_{+}}
\def\ep{\varepsilon}
\def\vp{\varphi}

\def\mC{\mathcal{C}}
\def\mL{\mathcal{L}}
\def\vE{\mathcal{E}}
\def\mI{\mathcal{I}}
\def\mJ{\mathcal{J}}
\def\supp{\text{supp}}

\def\mL{\mathcal{L}}
\DeclareMathOperator*{\essup}{\text{ess}\,\text{sup}}

\begin{document}

	
\title{Asymptotic stability of 3D out-flowing compressible viscous fluid under non-spherical perturbation}

\date{}
\author[$1$]{Yucong {\sc Huang}\footnote{Corresponding Author. Email Adress: huang@c.titech.ac.jp}}

\author[$1$]{Shinya {\sc Nishibata}
}

\affil[$1$]{\small Department of Mathematical and Computing Sciences, \protect\\ Institute of Science Tokyo, Tokyo 152-8552, Japan}

\maketitle

	\vspace{-6mm}
	
\begin{abstract}
    We study an outflow problem for the $3$-dimensional isentropic compressible Navier-Stokes equations. The fluid under consideration occupies the exterior domain of the unit ball $\Omega=\{x\in\R^3\,\vert\, |x|\ge 1\}$ and it is flowing out from the unit ball
    $\Omega$ at a constant speed $|u_b|$, in the normal direction to the boundary surface $\d\Omega$. The existence of a unique spherically symmetric stationary solution $(\tilde{\rho},\tilde{u})$ is obtained by I.~Hashimoto and A.~Matsumura in 2021, provided that the fluid velocity at the far-field is assumed to be zero, and $|u_b|$ is sufficiently small. Subsequently, authors of the present article prove in 2024 that $(\tilde{\rho},\tilde{u})$ is time-asymptotically stable under large spherically symmetric initial perturbations in the suitable Sobolev norm. The main purpose of the present paper is to investigate the case when the initial perturbations are possibly non-spherically symmetric. We show that $(\tilde{\rho},\tilde{u})$ remains asymptotically stable in time, under general small initial perturbations in the $H^3$-norm.
\end{abstract}

\vspace{-6mm}
	
	{\hypersetup{linkcolor=blue}
		\tableofcontents
	}

\section{Introduction}\setcounter{equation}{0}
\paragraph{} The Navier-Stokes equations for the isentropic motion of compressible viscous gas in the Eulerian coordinate is the system of equations given by
\begin{subequations}\label{ns}
\begin{gather}
\d_t \rho + \dd (\rho u) = 0, \label{ns1} \\
\d_t( \rho u) + \dd ( \rho u \otimes u) + \nabla P(\rho) = \mu \Delta u +(\mu + \lambda) \nabla \dd u .\label{ns2}
\end{gather}
\end{subequations}
We study the asymptotic behavior of a solution $(\rho, u)$ to (\ref{ns}) in an unbounded exterior domain $\Omega \vcentcolon= \{ x=(x_1,x_2,x_3)^\top \in \R^3 \; \vert \; |x| > 1 \}$. Here $\rho>0$ is the mass density; $u=(u^1,u^2,u^3)^{\top}$ is the velocity of gas; $P(\rho) = K\rho^\gamma \; (K>0, \gamma \ge 1)$ is the pressure with the adiabatic exponent $\gamma$; $\mu$ and $\lambda$ are constants called viscosity-coefficients satisfying $\mu>0$ and $2\mu + 3 \lambda \ge 0$. Hereafter, we denote $Lu = \big((Lu)^1,(Lu)^2,(Lu)^3\big)^{\top}$ with
\begin{align*}
    (Lu)^i \vcentcolon =  \mu \Delta u^i +(\mu + \lambda) \d_{x_i} \dd u =  \sum_{j=1}^3  \{ \mu \d_{x_j}^2 u^i + (\mu + \lambda) \d_{x_i}  \d_{x_j} u^j\} \ \text{ for } \ i=1,2,3. 
\end{align*}

In the present paper
we study the 
problem where fluid  flows out at a constant velocity from an inner sphere
of which 
center is at the origin with radius $|x|=1$. The corresponding boundary condition is given as
\begin{subequations}\label{bdry}
\begin{align}
&\ \, u(t,x) = u_b \frac{x}{|x|} \ \ &&\text{ for } \ x\in\d\Omega, \  &&\text{ where } \ u_b<0,\label{ub}\\
&\lim\limits_{|x|\to \infty} (\rho,u)(t,x) = (\br,0) \ \ &&\text{ for } \ t\ge 0, \  &&\text{ where } \ \br>0. \label{ff}
\end{align}  
\end{subequations}
For equations (\ref{ns})--(\ref{bdry}), we prescribe the initial data $(\rho_0,u_0)$ as
\begin{subequations}\label{iac}
	\begin{gather}
		\rho(x,0) = \rho_0(x)>0, \quad u(x,0) = u_0(x), \label{ic}\\
		\lim_{|x| \to \infty} (\rho_0(x), u_0(x)) = (\br, 0), \quad 0<\br <\infty. \label{asymp}
	\end{gather}
\end{subequations}
In addition, 
the data
(\ref{ic}) is assumed
to be
compatible with the boundary data (\ref{ub}).
Namely, we assume that
\begin{subequations}\label{compa}
\begin{gather} 
u_0\vert_{x\in \d\Omega} = u_b \frac{x}{|x|}\big\vert_{x\in \d\Omega} , \label{compa1} \\ 
\big\{-\rho_0 (u_0\cdot \nabla) u_0 + \mu \Delta u_0 + (\mu+\lambda)\nabla \dd u_0 - \nabla P(\rho_0) \big\} \big\vert_{x\in \d\Omega} = 0. \label{compa2}
\end{gather}
\end{subequations}
Since the characteristic velocity of (\ref{ns1}) is directed towards the origin at the boundary due to (\ref{ub}), no boundary condition for $\rho$ at $x\in\d\Omega$ is necessary for the well-posedness of the initial boundary value problem (\ref{ns})--(\ref{compa}).
\subsection{Spehrically symmetric stationary solution} 
\paragraph{} In this subsection, we give an overview of the spherically symmetric stationary solution to the problem \eqref{ns}--\eqref{bdry}. Consider the general multi-dimensional exterior domain $\Omega=\{x\in\R^n \,\vert\, |x|\ge 1\}$ where $n\ge 2$. Letting $(\tr,\tu)$ be the stationary solution of \eqref{ns}--\eqref{compa}. Solving the system
\begin{subequations}\label{st}
    \begin{align}
        &\dd \big( \tr \tu \big) = 0 && \text{for } \ x\in\Omega,\label{st1}\\
        &\tr (\tu\cdot \nabla ) \tu + \nabla P(\tr) = L \tu && \text{for } \ x\in\Omega,\label{st2}\\
        & \tu(x) = u_b \frac{x}{|x|} && \text{for } \ x\in\d\Omega,\label{stbdry}\\
        &(\tr,\tu)(x) \to (\br,0) && \text{as } \ |x|\to\infty. \label{stff}
    \end{align}
\end{subequations}
Since the boundary and far-field conditions \eqref{stbdry}--\eqref{stff} are spherically symmetric, if the unique solution $(\tr,\tu)$ to \eqref{st} exists, then $(\tr,\tu)$ must also be spherically symmetric. This means that $(\tr,\tu)$ is expressed in the form:
\begin{equation}\label{sph-st}
    \tr(x) = \tr(r), \quad \tu(x) = \tU(r) \frac{x}{r}, \qquad \text{where } \ r\vcentcolon = |x|\ge 1,
\end{equation}
for a certain function
$\tU\vcentcolon[1,\infty)\to \R$. Then the system \eqref{st} is reformulated as 
\begin{subequations}\label{rst}
\begin{align}
&\d_r \big(r^{n-1} \tr \tU \big) = 0 && \text{for } \ r\ge 1, \label{rst1} \\
&\tr \tU \d_r\tU + \d_r P(\tr) = (2\mu+\lambda) \d_r\Big( \frac{(r^{n-1} \tU)_r}{r^{n-1}} \Big) && \text{for } \ r\ge 1, \label{rst2}\\
&\tU(1) = u_b, \quad \lim_{r \to \infty} \big(\tr(r), \tU(r)\big) = (\br, 0). &&\label{rstbdry}
\end{align}
\end{subequations}
Multiplying the equation (\ref{st1}) by $r^{n-1}$, and integrating the resultant equality over $(1,r)$, we obtain
\begin{equation*}
\tU(r) = \dfrac{u_b \tr(1) }{r^{n-1}\tr(r)},
\end{equation*}
for $r\ge1$. If a solution $(\tr,\tU)(r)$ exists and $\tr(r)\to\br \in (0,\infty)$ as $r\to\infty$, then it is necessary to have the far-field condition, $\tU(r)\to 0$ as $r\to\infty$.

For general spatial dimension $n\ge 2$, I.~Hashimoto and A.~Matsumura in \cite{h-m21} proved the existence and uniqueness of a spherically symmetric classical solution to \eqref{st}, provided that $|u_b|$ is small. Subsequently in \cite{h-n-s23}, I.~Hashimoto, S.~Nishibata and S.~Sugizaki refined the spatial decay estimate for $(\tr,\tU)$, and they also prove the asymptotic stability of $(\tr,\tU)$ for small, spherically symmetric initial perturbations. Specifically, the well-posedness and properties of the stationary solution is summarised as follows:
\begin{lemma}\label{lem:st}
There exists $\delta=\delta(\br,\mu,\gamma,K,n)>0$ such that if $| u_b | \le \delta$, then the unique solution $(\tr,\tU)(r)\in \mC^2[1,\infty)$ to the problem (\ref{rst}) exists. In addition, $\tr(r)$ and $\tU(r)$ are respectively strictly monotone increasing and decreasing. Moreover, there exists a constant $C=C(\br,\mu,\gamma,K,n)>0$ such that for $r\ge1$,
\begin{equation*}
\begin{aligned}
&\tU(r) = u_b \tr(1) \dfrac{r^{1-n}}{\tilde{\rho}(r)}, \quad && | \tr(r) - \br | \le C r^{-2n+2},\\
&C^{-1} |u_b| r^{-n} \le |\d_r\tU(r)| \le C |u_b| r^{-n}, \quad && C^{-1} |u_b|^2 r^{-2n+1} \le |\d_r\tr(r)| \le C |u_b|^2 r^{-2n+1},\\
&C^{-1} |u_b| r^{-n-1} \le |\d_r^2\tU(r)| \le C |u_b| r^{-n-1}, \quad && C^{-1} |u_b|^2 r^{-2n} \le |\d_r^2\tr(r)| \le C |u_b|^2 r^{-2n}.
\end{aligned}
\end{equation*}
\end{lemma}
Y. Huang and S. Nishibata in \cite{HN-out} prove that if the initial data \eqref{iac} is spherically symmetric, then the stationary solution $(\tr,\tu)$ obtained in Lemma \ref{lem:st} asymptotically stable even for the 
large
initial perturbation. However, the stability under general possibly non-spherically symmetric initial perturbation 
has 
remained an open question. In the present paper, we give a positive answer to this problem.
Precisely
 we prove that the asymptotic stability holds 
provided that the initial perturbation $(\rho_0-\tr,u_0-\tu)$ is small in the Sobolev space, $H^3(\Omega)$.

\subsection{Main theorem}
\begin{theorem}\label{thm:main}
    There exists $\ep=\ep(\br,\mu,\gamma,K,n)>0$ such that if 
    \begin{equation*}
        |u_b| + \|(\rho_0-\tr,u_0-\tu)\|_{H^3(\Omega)} \le \ep,  
    \end{equation*}
    then the initial boundary value problem \eqref{ns}--\eqref{compa} has a unique solution $(\rho,u)$ globally-in-time
    in the function space
    \begin{gather*}
        \rho -\tr \in \mC^0\big([0,\infty);H^3(\Omega)\big)\cap \mC^1\big([0,\infty);H^2(\Omega)\big)\cap L^2\big(0,\infty;H^3(\Omega)\big),\\
        u-\tu \in \mC^0\big([0,\infty);H^3(\Omega)\big)\cap \mC^1\big([0,\infty);H^1(\Omega)\big)\cap L^2\big(0,\infty;H^4(\Omega)\big). 
    \end{gather*}
    In addition, the solution $(\rho,u)$ converges to
    the stationary solution
    $(\tr, \tu)$ as time goes to infinity, namely
    \begin{equation}\label{tacvg}
        \lim\limits_{t\to\infty} \sup\limits_{x\in\Omega} |(\rho-\tr,u-\tu)(t,x)| = 0.
    \end{equation}
\end{theorem}

\paragraph{Related results.}
The well-posedness problem of compressible Navier-Stokes equation with out-flowing or in-flowing boundary conditions have gained traction in recent decades. We state several previous results which are relevant to the present paper.

First, we refer readers to the book \cite{kaz90} by S.~N. Antontsev, A.~V. Kazhikhov, and V.~N. Monakhov, where a comprehensive survey of the mathematical theory of compressible Navier-Stokes equation is given. For the heat-conducing model, A.~Tani in \cite{Tani} provides a comprehensive proof for the existence and uniqueness of a local-in-time classical solution in the H\"older sapce. The major breakthrough research on the time-asymptotic stability of solution is first accomplished by A.~Matsumura and T.~Nishida in \cite{m-n83}, where they consider the heat-conductive compressible flow for a general $3$-dimensional exterior domain with external force and adhesion boundary condition ($u\vert_{\d\Omega}=0$ and $\d_{n}\theta\vert_{\d\Omega}=0$ where $\theta$ is the absolute temperature). In this case, the stationary solution takes the form $(\tilde{\rho},0,\bar{\theta})$ where $\tilde{\rho}=\tilde{\rho}(x)$ is a positive function of spatial variable, and $\bar{\theta}$ is a positive constant. Under smallness assumptions on the initial perturbation, it is proved that $(\tilde{\rho},0,\bar{\theta})$ is stable globally in time.

For spherically symmetric flow in the bounded annulus domains with adhesion boundary condition, a pioneering work has been done by N.~Itaya \cite{itaya85}, which establishes the global-in-time existence of a unique classical solution for the isothermal model. In this work, no smallness assumption on the initial data is imposed. Later, T.~Nagasawa investigates the asymptotic stability for the same problem in \cite{nagasawa}. The paper \cite{itaya85} has motivated a sequence of developments on the topic of spherically symmetric solution. For instance, A.~Matsumura in \cite{matsu92} constructs a spherically symmetric classical solution to the isothermal model with external forces on a bounded annular domain. Moreover, he also proves that the convergence rate for the time-asymptotic stability is exponential. Subsequently, these result is extended to the isentropic and heat-conductive models by K.~Higuchi in \cite{higuchi92}. The well-posedness of a spherically symmetric solution in an unbounded exterior domain is first obtained by S.~Jiang \cite{jiang96}, where the global-in-time existence of the uniqueness classical solution is shown. In addition, a partial result on the time asymptotic stability is proved in \cite{jiang96} where, for $n=3$, $\|u(\cdot,t)\|_{L^{2j}}\to 0$ as $t\to\infty$ with an arbitrary fixed integer $j\ge 2$. Later, this restriction on the long time stability was fully resolved by T.~Nakamura and S.~Nishibata 
in \cite{n-n08}, where a complete stability theorem is obtained for large initial data. We also refer to the paper \cite{NNY} by T.~Nakamura, S.~Nishibata and S.~Yanagi, where the time asymptotic stability of the spherically symmetric solution for the isentropic flow is established with large initial data and external forces.

For the one-dimensional outflow or inflow problem, the stationary solution becomes non-trivial. This leads to a variety of physically interesting time-asymptotic behaviours for the solutions. A.~Matsumura in \cite{m2001} starts the first investigation of this problem for the isentropic model posed on the one dimensional half-space domain. Several kinds of boundary conditions are considered in \cite{m2001}, which includes inflow, outflow and adhesion boundary conditions. He formulates
conjectures on the classification of 
asymptotic behaviours of the solutions in different cases subject to the relation between the boundary data and the spatial asymptotic data. Then the stability theorems for several cases of inflow problem were established by A.~Matsumura and K. Nishihara in \cite{m-n01},
where they employ the Lagrangian mass coordinate.
Following this work, S.~Kawashima, S.~Nishibata and P.~Zhu in \cite{k-n-z03} further study the outflow problem in one dimensional half space. They prove the long time stability of solutions with small initial perturbation 
from the stationary solution. 
A detailed desicription of the convergence rate towards the stationary
solution is found in the work of T.~Nakamura, S.~Nishibata and T.~Yuge \cite{n-n-y07}. For the non-isentropic inflow problem in the half line, T.~Nakamura and S.~Nishibata in \cite{NN} demonstrate the time-asymptotic stability of stationary solutions under a small initial perturbation, for both the subsonic and transonic cases. 

The research on the outflow and inflow problems for fluid flow in multi-dimensional (M-D for short) domains are relatively new subjects. For isentropic fluid occupying the M-D half space $\mathbb{R}^n_{+}$, Y.~Kagei and S.~Kawashima \cite{k-k2006} show that the planar stationary solution for the out-flowing boundary condition is asymptotically stable. Subsequently T.~Nakamura and S.~Nishibata in \cite{n-n09} obtain the convergence rates for this time-asymptotic stability. In particular, they show that an algebraic or an exponential decay rate holds for the supersonic flow. Precisely,
an algebraic convergence rate holds for the transonic flow
provided that the respective initial perturbations decay in a normal direction to the boundary surface with the corresponding decay rates. For isentropic fluid occupying the exterior domain of a unit ball, I.~Hashimoto and A.~Matsumura in \cite{h-m21} obtain the existence of the unique spherically symmetric stationary solutions for both inflow and outflow problems provided that the velocity at boundary, $|u_b|$ is sufficiently small. More recently, I.~Hashimoto, Y.~Huang and S.~Nishibata in \cite{h-h-n24} study the spherically symmetric inflow problem in the exterior domain and they prove that if the far-field state of the fluid is not vacuum but suitably rarefied, then the corresponding stationary solution is stable under small initial perturbation. For the exterior outflow problem, I.~Hashimoto, S.~Nishibata and S.~Sugizaki in \cite{h-n-s23} show the stability of stationary solution with small spherically symmetric initial perturbations. Subsequently, this smallness assumption on initial data is removed by Y. Huang and S. Nishibata in \cite{HN-out} and they prove that the stationary solution for the outflow problem remains asymptotically stable under large spherically symmetric initial perturbations. In addition, E. Feireisl, P. Gwiazda, and A. Swierczewska-Gwiazda consider the finite energy weak solution to the exterior outflow problem in \cite{FGG}. The main theorem proved in the present paper is the asymptotic stability for the outflow problem under general possibly non-spherically symmetric initial perturbations in
the exterior domain.

\paragraph{Outline of the paper.}
First, in Section \ref{ssec:reform}, we restate equations \eqref{ns} in terms of the difference functions $(\phi,\psi)=(\rho-\tr,u-\tu)$ and we also reformulate the equations in linear forms which are suitable for the a-priori estimates. Moreover, we provide the precise definition of energy norm $N(T)$ which is used throughout the proof. Subsequently, in Section \ref{ssec:RE}, we derive the $L^2$--estimates of $(\phi,\psi)$ by constructing the $L^2$--relative energy density, which is denoted by $\vE$. Sections \ref{ssec:ETD}--\ref{ssec:stokes} are devoted to the derivation of $L^{\infty}_{T}H^2$-estimate for $\big(\nabla\phi,\nabla \psi\big)$. Note that out-flowing boundary condition \eqref{bdry} plays a significant role in the derivation for these estimates involving higher-order derivatives. Since $\d\Omega$ is a sphere, we observe that the temporal and angular derivatives preserve the boundary condition \eqref{bdry}.
On the other hand, the radial derivative of $u$ at $\d\Omega$ is not a-priori obtained. Hence,
it is necessary to split the analysis of $L^{\infty}_{T}H^2$-estimate for $(\nabla\phi,\nabla\psi)$ into 3 steps as follows:\\
1. derive estimates for temporal and angular derivatives of $(\phi,\psi)$ by using the parabolic structure of \eqref{ns2}, which are the main contents of Sections \ref{ssec:ETD}--\ref{ssec:angular};\\
2. derive estimates for radial derivative of $\phi$ using the damping effect of compressible Navier-Stokes equations, which is illustrated in Section \ref{ssec:damp};\\
3. derive estimates for the full spatial derivatives of $(\phi,\psi)$ by treating the system \eqref{ns} as Stokes equations and this is devoted to Section \ref{ssec:stokes}. It is worth noting that in step 2 (Section \ref{ssec:damp}), a cancellation occurs on the radial derivatives of $\psi$. This serves as the key argument to successfully close the estimate. 

\paragraph{Notations.}
\begin{enumerate}
    \item For $T\ge 0$, denote the space-time domains as 
    \begin{equation*}
        Q_T \vcentcolon= [0,T]\times\Omega, \qquad \Gamma_T \vcentcolon= [0,T] \times \d\Omega.
    \end{equation*}
    \item For the Cartesian coordinate $(x_1,x_2,x_3)\in\R^3$, denote $\he_1$, $\he_2$, $\he_3$ as the corresponding unit basis vectors. Moreover, for functions $F\vcentcolon\R^3 \to \R$, and $V\vcentcolon \R^3 \to \R^3$ with $V=(V^1,V^2,V^3)^{\top}$, denote the differential operators:
    \begin{align*}
        &\nabla F \vcentcolon= \sum_{i=1}^3 \he_i \d_{x_i}F, \quad && \dd V\vcentcolon= \sum_{i=1}^3\d_{x_i} V^i, \quad && \Delta F \vcentcolon = \sum_{i=1}^3 \d_{x_i}^2 F,\\
        &(V\cdot \nabla )F \vcentcolon= \sum_{i=1}^3 V^i \d_{x_i} F, \quad && D_t F \vcentcolon= \d_t F + (u\cdot \nabla) F. &&
    \end{align*}
    In addition, for operators $S$, $T$, denote $[S,T]\vcentcolon= S T - T S$ as the commutator.
    \item For $k>0$, the Sobolev space is defined by
    \begin{align*}
        H^k(\Omega) \!=\! \bigg\{ F \in L^2(\Omega) \, \bigg\vert \, D^k F \vcentcolon= \Big\{ \frac{\d^\alpha F}{\d x_1^{\alpha_1}\d x_2^{\alpha_2}\d x_3^{\alpha_3}}, \, |\alpha|=l \Big\} \in L^2(\Omega), \ 1\le l \le k \bigg\},
    \end{align*}
    Furthermore, the norms are defined by
    \begin{align*}
        \|F\|_{k} \vcentcolon= \|F\|_{H^k(\Omega)}, \qquad \|F\|\vcentcolon= \|F\|_0 = \|F\|_{L^2(\Omega)}, \qquad \|F\|_{\infty} \vcentcolon= \essup\limits_{x\in\Omega} |F|.
    \end{align*}
\end{enumerate} 

\section{A-priori estimates}\label{sec:a-priori}
\setcounter{equation}{0}

\subsection{Reformulation and energy norm \texorpdfstring{$N(T)$}{N(T)}}\label{ssec:reform}
\paragraph{} In this subsection, we reformulate the problem \eqref{ns}-\eqref{iac} to a form which is suitable for the a-priori estimate. In addition, we also define the notion of energy norm $N(T)$.  

Define the difference functions
between the solution $(\rho, u)$ 
to \eqref{ns}-\eqref{iac} and the stationary solution $(\tr, \tu)$ to \eqref{rst} by
\begin{equation*}
    \phi\vcentcolon= \rho-\tr, \qquad \psi\vcentcolon= u-\tu, \qquad (\phi_0,\psi_0)\vcentcolon= (\phi,\psi)\vert_{t=0}
    =(\rho_0-\tr, u_0-\tu).
\end{equation*}
Then the system \eqref{ns}-\eqref{iac} is rewritten to
that of the unknown functions $(\phi, \psi)$
as 
\begin{subequations}\label{LE}
    \begin{align}
        \psi\vert_{\d\Omega} =& 0 \qquad \text{ for } \ t\ge 0,\label{DiffBdry}\\
        \mL^0(\phi,\psi)\vcentcolon=& D_t \phi + \br \dd \psi  = f^0,\label{E0}\\
        \mL^i(\phi,\psi)\vcentcolon=& \d_t \psi^i - \frac{\mu}{\br} \Delta \psi^i - \frac{\mu+\lambda}{\br} \d_{x_i}\dd \psi + \bq \d_{x_i}\phi = f^i,\label{Ei}
    \end{align}
    for $i=1,2,3$, where
    $D_t\vcentcolon= \d_t + u\cdot \nabla$ denotes the material derivative,
    $q(\rho)\vcentcolon= P^{\prime}(\rho)/\rho $,
    $\bq\vcentcolon=q(\br)$. Moreover, 
    \begin{align}
        f^0 \vcentcolon=& - \phi \dd \psi + (\br-\tr)\dd\psi - \psi\cdot \nabla \tr - \phi \dd \tu,\label{f0}\\
        f^{i} \vcentcolon=& -\frac{\phi+\tr-\br}{\br(\phi+\tr)} \big\{ \mu \Delta \psi^i + (\mu+\lambda) \d_{x_i} \dd \psi \big\} + \tilde{f}^i,\label{fi}\\ 
        \tf^{i} \vcentcolon=& -(\psi\cdot \nabla) \psi^i - (\tu\cdot \nabla) \psi^i - (\psi\cdot \nabla) \tu^i - \frac{\phi}{\phi+\tr} (\tu\cdot \nabla) \tu^i \nonumber\\
        &+ \{ \bq - q(\rho) \} \d_{x_i} \phi - \frac{P^{\prime}(\br)-P^{\prime}(\rho)}{\phi+\tr} \d_{x_i} \tr . \label{tfi}
    \end{align}
    In addition, we define $\mL := \hat{e}_1 \mL^1 + \hat{e}_2 \mL^2 + \hat{e}_3 \mL^3$ as the vector-valued operator and $f := \hat{e}_1 f^1 + \hat{e}_2 f^2 + \hat{e}_3 f^3$.
\end{subequations}
Moreover, for given initial data $(\rho_0,u_0)$ and time $T>0$, we define the solution space $X(T)$ as
\begin{equation*}
    X(T) \!\vcentcolon=\! \bigg\{ (\phi,\psi)\!\in\! L^{\infty}\big(0,T; H^3(\Omega)\big) \,\bigg\vert\, \begin{aligned}
        &\phi\!\in\! H^1\big(0,T;H^2(\Omega)\big),\,\psi \!\in\! H^1\big(0,T;H^1(\Omega)\big),\\
        &\psi\vert_{\d\Omega}\!=\!0,\, (\phi,\psi)\vert_{t=0}\!=\!(\phi_0,\psi_0),\, N(T)<\infty.  
    \end{aligned}\bigg\}, 
\end{equation*}
where the energy norm $N(T)$ is defined by
\begin{align}
    N(T) \vcentcolon=& \sup\limits_{0\le t \le T} \Big\{  \|(\phi,\psi)(t,\cdot)\|_{3} + \|\d_t \phi(t,\cdot)\|_{2} + \| \d_t \psi(t,\cdot) \|_{1}  \Big\} \nonumber\\
    &+ \bigg\{ \int_{0}^T\! \big\{ \| \nabla \phi \|_2^2 + \|\d_t \phi\|_{3}^2 +\|\nabla \psi\|_{3}^2 + \|\d_t \psi \|_{2}^2 \big\}\, \dif t\bigg\}^{\frac{1}{2}}.\label{norm}
\end{align}
To obtain the global-in-time well-posedness and time-asymptotic stability, we need the existence and uniqueness of a local-in-time solution. This is stated as follows.
\begin{lemma}\label{lem:local}
    Suppose the stationary solution $(\tr,\tu)$ and initial data $(\rho_0,u_0)$ satisfies $N(0) <\infty$. Then there exists $T_{\ast}>0$, depending only on $N(0)$, such that the initial-boundary value problem \eqref{ns}--\eqref{compa} has a unique strong solution $(\rho,u)$ satisfying $N(T)<\infty$ for $T\in[0,T_{\ast})$.
\end{lemma}
The proof of Lemma \ref{lem:local} follows the same procedure in \cite{m-n82}, and so we omit the details. We mention that for the Cauchy problem, an analogous method is employed in \cite{MN1980} to obtain the local-in-time well-posedness(See Theorem 5.2 of \cite{MN1980}).

For initial data with $N_0\equiv N(0)<\infty$, there exists a time $T=T(N_0)>0$ such that
\begin{equation} \label{continuity}
N(t)\le 2 N(0) \qquad \text{for } \ t\in [0,T]
\end{equation}
by the continuity of $t\mapsto N(t)$. By the Sobolev embedding theorem $H^2(\Omega) \xhookrightarrow[]{} \mC^0(\Omega)$, there exists a constant $\ep= \ep(\br,u_b,\gamma,K,\mu,\lambda)>0$ such that if $|u_b| + N(0) \le \ep$ then
\begin{equation}\label{priorirho}
    \frac{\br}{2} \le \rho(t,x) \le \frac{3}{2}\br \qquad \text{for } \ (t,x)\in [0,T]\times \Omega. 
\end{equation}
Moreover, it also follows by the Sobolev embedding theorem that there exists a constant $C=C(\br,u_b,\gamma,K,\mu,\lambda)>0$ such that,
\begin{subequations}\label{Linf}
    \begin{gather}
        \sup\limits_{0\le t\le T}\big\| \big( \phi, \nabla \phi, \psi, \nabla \psi, \d_t\phi \big) (t,\cdot) \big\|_{\infty} \le C N(T),\\
        \int_{0}^T \big\|\big(\nabla \phi,\d_t\phi , \nabla \d_t \phi , \nabla^2 \psi, \d_t \psi\big)(t,\cdot)\big\|_{\infty} \, \dif t \le C N(T).
    \end{gather}    
\end{subequations}
The main aim of 
Section \ref{sec:a-priori} is to derive the following a-priori estimate:
\begin{theorem}\label{thm:apriori}
    There exists two positive constants $\ep_0=\ep_0(\rho_+,\gamma,K,\mu,\lambda)>0$ and $C=C(\rho_+,\gamma,K,\mu,\lambda)>0$, which are independent of $|u_b|$, such that the following statement holds: if $(\rho,u)$ solves the problem \eqref{ns}--\eqref{compa} in some time interval $t\in [0,T]$ and its energy norm and boundary data satisfy $|u_b|+N(T)\le \ep_0$, then $N(T) \le C N(0).$
\end{theorem}
If Lemma \ref{lem:local} and Theorem \ref{thm:apriori} are obtained, then the global existence of $(\rho,u)$ is proved by the standard continuity argument. For the detail, we refer readers to Section 3 of \cite{m-n83}. Moreover, the asymptotic convergence \eqref{tacvg} is derived by the same argument presented in Theorem 7.1 of \cite{MN1980}.   

\subsection{\texorpdfstring{$L^2$}{L2}-estimate}\label{ssec:RE}
\paragraph{} To show the a-priori estimate,
we first derive the $L^2$-estimate
of the perturbation $(\phi,\psi)$.
To this end, we employ a potential energy density $H(\zeta,\xi)$,
which is defined by
\begin{equation}\label{H}
H(\zeta,\xi)\! \vcentcolon=\! \zeta\! \int_{\xi}^{\zeta}\! \dfrac{P(z)-P(\xi)}{z^2} \, \dif z\! =\! \left\{ 
\begin{aligned}
&\dfrac{K}{\gamma-1}\big\{ \zeta^{\gamma}-\xi^{\gamma} - \gamma \xi^{\gamma-1}  (\zeta-\xi) \big\} && \text{if } \ \gamma>1,\\
&K \xi \big( 1 - \dfrac{\zeta}{\xi} + \dfrac{\zeta}{\xi} \ln \dfrac{\zeta}{\xi} \big) && \text{if } \ \gamma=1,
\end{aligned}
\right.
\end{equation} 
where 
 $\zeta$, $\xi\ge 0$.
It is holds that $H(\zeta,\xi)$ is a positive convex function achieving global minimum at $\{\zeta=\xi\}$ with $H(\zeta,\xi)=0$ if and only if $\zeta=\xi$. Moreover $H(\zeta,\xi)$ satisfies  
\begin{subequations}\label{Hiden}
\begin{gather}
\zeta \d_\zeta H(\zeta,\xi) = H(\zeta,\xi) + P(\zeta) - P(\xi), \qquad \xi \d_{\xi} H(\zeta,\xi) = - P^{\prime}(\xi) (\zeta-\xi),\label{Hiden-1}\\
\zeta \d_\zeta H(\zeta,\xi) + \xi \d_{\xi} H(\zeta,\xi) =
\gamma H(\zeta,\xi), \label{Hiden-2}\\
P(\zeta)- P(\xi) - P^{\prime}(\xi)(\zeta-\xi) = (\gamma-1) H(\zeta,\xi).\label{Hiden-3}
\end{gather}
\end{subequations}
A $L^2$--relative energy density $\vE$ is defined by
\begin{equation}\label{RE}
\vE \vcentcolon= \dfrac{1}{2} \rho |u-\tu|^2 + H(\rho,\tr), \qquad \vE_0 := \vE\vert_{t=0}.
\end{equation}
By Taylor's theorem, we
have the following proposition.
\begin{proposition}\label{prop:ee}
If there exist constants $0<\rho_{\ast} < \rho^{\ast}<\infty$ such that $\rho_{\ast} \le \rho,\, \tr  \le \rho^{\ast}$, then 
\begin{equation*}
C^{-1}|\rho-\tr|^2 \le H(\rho,\tr) \le C |\rho-\tr|^2,
\end{equation*}
for a certain constant $C=C(\rho^{\ast},\rho_{\ast},\gamma,K)>0$.
\end{proposition}
\begin{lemma}\label{lemma:rE}
There exists $\delta_0=\delta_0(\br,\gamma,K,\mu,\lambda)>0$ such that if $|u_b|\le \delta_0$ then, 
\begin{align*}
&\sup\limits_{0\le t\le T} \int_{\Omega} \vE(t,x)\,\dif x + \iint_{Q_T} \big\{ \dfrac{\mu}{2} |\nabla\psi|^2 + (\mu+\lambda)|\textnormal{\dd} \psi|^2 \big\}\,\dif x \dif t && \nonumber\\	
&+ |u_b|\int_{0}^T\!\!\!\!\int_{\{|y|=1\}}\!\!\!\!\!H(\rho,\tr)\,\dif\mathcal{S}_{y}\dif t + |u_b|^3 \iint_{Q_T}  \dfrac{|\phi|^2}{|x|^{7}}\,\dif x \dif t + |u_b|\iint_{Q_T}\!\!  \frac{|x\cdot \psi|^2}{|x|^9}  \, \dif x \dif t \le C \vE_0,
\end{align*}
for $T\in[0,T_0]$, where $T_0$ is given in \eqref{continuity}--\eqref{priorirho}, and $C=C(\br,\gamma,K,\mu,\lambda)>0$ is a certain constant independent of $u_b$ and $(\rho_0,u_0)$.
\end{lemma}
\begin{proof}
Equations \eqref{ns1} and \eqref{st1} yields that
\begin{align*}
&\d_t H(\rho,\tr) + \big( P(\rho)-P(\tr) \big) \dd \psi\\ 
=&(1-\gamma) H(\rho,\tr) \dd \tu  -\dd \big(u H(\rho,\tr)\big) - \dfrac{\phi}{\tr} (\psi\cdot \nabla) P(\tr).
\end{align*}
Integrating the above 
equality
over $\Omega$ with
using the boundary condition \eqref{ub}, we obtain
\begin{align}\label{sHtInt}
    &\dfrac{\dif }{\dif t} \int_{\Omega}\! H(\rho,\tr)(t,x)\, \dif x + \int_{\Omega}\!\! \big( P(\rho)-P(\tr) \big) \dd \Psi \dif x - u_b \int_{\{|y|=1\}}\!\!\!\! H(\rho,\tr)(t,y)\,\dif S_{y} \nonumber\\
    &= (1-\gamma) \int_{\Omega}\!\! H(\rho,\tr)\dd \tu\,\dif x - \int_{\Omega} \dfrac{\phi}{\tr} (\psi\cdot\nabla) P(\tr) \,\dif x. 
\end{align}
Using \eqref{ns} and \eqref{st}, we have
\begin{align*}
	&\d_t \big(\rho\dfrac{|\psi|^2}{2}\big)  - \big(P(\rho)-P(\tr)\big) \dd \psi + \mu|\nabla\psi|^2 + (\mu+\lambda)|\dd\psi|^2 + \dd \mathcal{F}\\
 =& -\psi \cdot \big\{ \phi (\psi\cdot\nabla) \tu + \phi (\tu\cdot\nabla) \tu + \tr(\psi\cdot\nabla) \tu\big\},
\end{align*}
where
\begin{align*}
    \mathcal{F}\vcentcolon= \rho u \dfrac{|\psi|^2}{2} + \psi \big( P(\rho) - P(\tr) \big) - \mu(\nabla \psi)^{\top} \cdot \psi - (\mu+\lambda) \psi \dd \psi.
\end{align*}
Integrating the above equality over $\Omega$ with using the boundary condition \eqref{bdry}, we have
\begin{align}\label{spsiInt}
	&\dfrac{\dif}{\dif t} \int_{\Omega}\!\! \rho \dfrac{|\psi|^2}{2}\,\dif x- \int_{\Omega}\!\! \big( P(\rho)-P(\tr) \big) \dd\psi\,\dif x + \mu\int_{\Omega}\!\! |\nabla\psi|^2\,\dif x  + (\mu+\lambda) \int_{\Omega}\!\! |\dd \psi|^2\,\dif x  \nonumber\\
	&=-\int_{\Omega}\!\! \psi\cdot \big\{ \phi (\psi\cdot\nabla) \tu + \phi (\tu\cdot\nabla) \tu + \tr(\psi\cdot\nabla) \tu\big\} \, \dif x. 
\end{align}
The integrands of the left hand of 
\eqref{spsiInt} is rewritten as
\begin{align*}
&\psi\cdot \big\{ \phi (\psi\cdot\nabla) \tu + \phi (\tu\cdot\nabla) \tu + \tr(\psi\cdot\nabla) \tu\big\} + \dfrac{\phi}{\tr} (\psi\cdot\nabla)P(\tr) = \rho \psi\cdot (\psi \cdot \nabla) \tilde{u} + \phi \psi \cdot \dfrac{L\tu}{\tr} .
\end{align*}
The above equality is obtained
from
the stationary momentum equation \eqref{st2}.
Adding equation (\ref{sHtInt}) to (\ref{spsiInt}) with using definition \ref{RE}, we obtain
that
\begin{align}\label{Et'}
	&\dfrac{\dif}{\dif t} \int_{\Omega}\! \vE\,\dif x + |u_b|\!\!\int_{\{|y|=1\}}\!\!\! H(\rho,\tr)(t,y)\,\dif\mathcal{S}_{y}+\int_{\Omega}\big\{ \mu |\nabla\psi|^2 + (\mu+\lambda)|\dd \psi|^2 \big\}\,\dif x \nonumber\\
	= & (1-\gamma) \int_{\Omega}\!\! H(\rho,\tr) \dd \tu\,\dif x - \int_{\Omega}\!\!\big\{ \rho   \psi^{\top} \cdot \nabla\tilde{u} \cdot \psi + \phi \psi \cdot \dfrac{L\tu}{\tr} \big\}\,\dif x.
\end{align}
The right hand side of \eqref{Et'}
is handled as follows.
As $\tu=\frac{x}{r}\tU(r)$ for $r=|x|$, $\nabla \tu$ is a symmetric matrix and
satisfies 
\begin{gather*}
\nabla \tu
=(\nabla \tu)_{+}+(\nabla \tu)_{-} \quad \text{and} \quad \dd \tu = \d_r \tU(r) + (n-1) \dfrac{\tU(r)}{r}, 
\end{gather*}
where
\begin{gather*}
    (\nabla \tu)_{+} \vcentcolon= \big( \d_r \tU(r) - \dfrac{\tU(r)}{r} \big) \dfrac{x\otimes x}{r^2}, \qquad  (\nabla \tu)_{-} \vcentcolon= \dfrac{\tU(r)}{r}\mathbb{I}_3.
\end{gather*}
By Lemma \ref{lem:st}, we have
the inequality
\begin{equation}
\dd \tu = \dfrac{\d_r (r^{2}\tU)}{r^{2}} = - \tr(1) u_b \dfrac{\d_r \tr(r)}{r^{2}|\tr(r)|^2} = \tr(1)|u_b| \dfrac{\d_r \tr(r)}{r^{2}|\tr(r)|^2}\ge \frac{|u_b|}{C r^7}.
\end{equation}
Hence using the above
inequalities
and Lemma \ref{lem:st}, we have that for an arbitrary vector $V\in\R^3$,
\begin{align*}
V^{\top}\cdot (\nabla \tu)_{+} \cdot V =&\big\{ \d_r \tU -\dfrac{\tU}{r} \big\} \big|\dfrac{x}{r}\cdot V\big|^2\\
=& \big\{ \dd \tu + \dfrac{3|u_b|\tr(1)}{r^3 \tr(r)} \big\}\big|\dfrac{x}{r}\cdot V\big|^2 \ge \frac{|u_b|}{C} \frac{|x\cdot V|^2}{r^9}.
\end{align*}
Substituting this
inequality
in (\ref{Et'}), we obtain that
\begin{align}\label{Et}
&\dfrac{\dif}{\dif t} \int_{\Omega} \vE\,\dif x + |u_b|\!\!\int_{\{|y|=1\}}\!\!\! H(\rho,\tr)(t,y)\,\dif\mathcal{S}_{y}+\int_{\Omega}\big\{ \mu |\nabla\psi|^2 + (\mu+\lambda)|\dd \psi|^2 \big\}\,\dif x \nonumber\\	
&+(\gamma-1) \int_{\Omega}\!\!  H(\rho,\tr)\dd \tu\,\dif x +\frac{|u_b|}{C} \int_{\Omega} \frac{|x\cdot \psi|^2}{|x|^9} \, \dif x\nonumber\\
\le& - \int_{\Omega}\!\! \rho \psi^{\top} \cdot (\nabla \tu)_{-}\cdot \psi\, \dif x  - \int_{\Omega}\!\!\phi \psi \cdot \dfrac{L\tu}{\tr}\,\dif x.
\end{align}
We estimate the right hand side of
\eqref{Et} term by term.
By Lemma \ref{lem:st}, \eqref{priorirho}, Propositions \ref{prop:ee} and \ref{prop:hardy}, the first term is estimated as
\begin{align}\label{temp1:re1}
    \Big| \int_{\Omega} \rho \psi^{\top}\cdot (\nabla \tu)_{-}\cdot \psi\,\dif x \Big|
    \le& \frac{3\br |u_b|}{2} \int_{\Omega} \dfrac{|\psi|^2}{|x|^3}\, \dif x \nonumber\\ 
    \le& 9 \br |u_b| \int_{\Omega}|\nabla \psi|^2 \, \dif x\le  \dfrac{\mu}{4}\int_{\Omega} |\nabla\psi|^2\, \dif x,
\end{align}
where we have chosen as
$|u_b|\le \mu/(36\br)$. In addition, It holds by Lemma 
\ref{lem:st}
that
\begin{align*}
    &\dd \tu = \tr(1)|u_b|\dfrac{\partial_r \tr}{r^{2}\tr^2} \ge \frac{|u_b|^3}{C r^7}, \quad 
    |L\tu| 
    = \mu \tr(1) |  u_b |\cdot \Big| \dfrac{\partial_r^2 \tr}{r^{2}\tr^2} -2 \dfrac{\partial_r \tr }{r^{3}\tr^2} - 2 \dfrac{|\partial_r \tr|^2}{r^{2}\tr^3} \Big| \le \frac{C|u_b|^3}{r^8}.
\end{align*}
By
using the above inequality and Proposition \ref{prop:ee}, we obtain that
\begin{align}\label{temp2:re1}
	(\gamma-1)\int_{\Omega} H(\rho,\tr) \dd \tu \,\dif x \ge \frac{|u_b|^3}{C_1} \int_{\Omega} \dfrac{|\phi|^2}{|x|^7}\,\dif x
\end{align}
for a certain constant $C_1=C_1(\br,\mu,\gamma,K)>0$. Moreover, owing to Lemma \ref{lem:st}, Cauchy-Schwarz's inequality and Proposition \ref{prop:hardy}, we get
\begin{align}\label{temp3:re1}
    &\Big|\int_{\Omega}\!\!\phi \psi \cdot \dfrac{L\tu}{\tr}\,\dif x \Big| \le C |u_b|^3 \int_{\Omega} \dfrac{|\phi \psi|}{|x|^{8}}\, \dif x  
    \le \dfrac{\mu}{4}\int_{\Omega}\!\! |\nabla\psi|^2\, \dif x + \frac{|u_b|^3}{C_1}\int_{\Omega} \dfrac{|\phi|^2}{|x|^{7}}\, \dif x,
\end{align}
where we 
have
chosen $|u_b|$ so small that $C|u_b|^3 \le 1/C_1$. Substituting (\ref{temp1:re1})--(\ref{temp3:re1}) in (\ref{Et}), we obtain
\begin{align*}
	&\dfrac{\dif}{\dif t} \int_{\Omega}\!\! \vE\,\dif x +\int_{\Omega}\!\!\big\{ \dfrac{\mu}{2} |\nabla\psi|^2 + (\mu+\lambda)|\dd \psi|^2 \big\}\,\dif x \nonumber\\	
	&+ |u_b|\int_{\{|y|=1\}}\!\!\! H(\rho,\tr)(t,y)\,\dif\mathcal{S}_{y} + \dfrac{|u_b|^3}{2C_1} \int_{\Omega}  \dfrac{|\phi|^2}{|x|^{7}}\,\dif x + \frac{|u_b|}{C} \int_{\Omega} \frac{|x\cdot \psi|^2}{|x|^9} \, \dif x \le 0.
\end{align*}
Integrating the above
inequality 
in time, we obtain the desired estimate.
\end{proof}

\subsection{Estiamtes for temporal derivatives}\label{ssec:ETD}
\paragraph{} In this section, we obtain estimates for the temporal derivatives of $(\phi,\psi)$. We remark that the temporal derivative preserves the boundary condition \eqref{DiffBdry} in the sense that
    \begin{equation}\label{dt-bdry}
        \d_t\psi(t,x) = \d_t^2 \psi(t,x) = 0 \qquad \text{for } \ (t,x)\in \Gamma_T \vcentcolon = [0,T]\times \d\Omega.
    \end{equation}
\begin{lemma}\label{lemma:timeD}
    Suppose $(\phi,\psi)\in X(T_0)$ is a solution to the system \eqref{LE}, where $T_0>0$ is the time given in \eqref{continuity}--\eqref{priorirho}. Then there exists a constant $C=C(\br,\gamma,K,\mu,\lambda)>0$ such that for $T\in[0,T_0]$ and $k=0,1$,
    \begin{align*}
        &\sup\limits_{0\le t\le T} \big\{ \| \d_t (\phi,\psi) \|^2 + \|\nabla \d_t^k \psi\|^2 \big\} +\!\! \int_{0}^T \!\!\!\! \big\{ \| \nabla \d_t \psi \|^2 \!+\! \| \d_t^{k+1}(\phi,\psi) \|^2 \big\}\dif t  + \!\! \iint_{\Gamma_T}\!\! |\d_t \phi|^2 \dif x \dif t \\
         \le& C \big\{ N^2(0) + |u_b| N^2(T) + N^3(T) \big\}. 
    \end{align*}
\end{lemma}
\begin{proof}
    Taking the derivative $\d_t$ on \eqref{E0} and multiplying the resultant equality by $\frac{\bq}{\br}\d_t\phi$, we get 
    \begin{align}
        &\d_t\big(\frac{\bq}{2\br}|\d_t\phi|^2\big) + \dd\big( \frac{\bq}{2\br} u |\d_t\phi|^2 \big) + \bq  \d_t\phi \dd \d_t\psi \nonumber \\
        =& \frac{\bq}{\br} \big\{ \d_t\phi \d_t f^0 + \frac{|\d_t\phi|^2}{2}\dd u  -\d_t\phi (\d_t\psi \cdot \nabla ) \phi\big\} . \label{temp:timeD1}
    \end{align}
    In addition, taking the derivative $\d_t$ on \eqref{Ei}, multiplying the resultant equality by $\d_t\psi$, and integrating by parts, we get
    \begin{align}
        &\d_t\big( \frac{1}{2}|\d_t\psi|^2 \big) + \frac{\mu}{\br} |\nabla \d_t\psi|^2 + \frac{\mu+\lambda}{\br} |\dd \d_t\psi|^2 - \bq \d_t\phi \dd\d_t\psi \nonumber \\
        =& \d_t\psi\cdot \d_t f + \dd \Big\{ \frac{\mu}{\br} \sum_{i=1}^3\d_t\psi^i \nabla \d_t\psi^i + \frac{\mu+\lambda}{\br} \d_t\psi \dd \d_t\psi  - \bq \d_t\psi \d_t\phi \Big\}. \label{temp:timeD2}
    \end{align}
    Summing equations \eqref{temp:timeD1}--\eqref{temp:timeD2}, integrating the resultant equality in $(t,x)\in Q_T$, applying the boundary condition \eqref{dt-bdry}, Lemma \ref{lem:st}, and inequality \eqref{Linf}, we obtain
    \begin{align}
        &\frac{1}{2} \int_{\Omega}\!\! \big\{ \frac{\bq}{\br}|\d_t\phi|^2 + |\d_t\psi|^2 \big\}\dif x \bigg\vert_{t=0}^{t=T} +  \frac{|u_b|}{2}\iint_{\Gamma_T} \!\! |\d_t\phi|^2 \, \dif x \dif t  \nonumber\\ 
        &+ \iint_{Q_T} \!\! \Big\{ \frac{\mu}{\br} |\nabla \d_t\psi |^2 + \frac{\mu+\lambda}{\br} |\dd \d_t\psi|^2 \Big\} \dif x \dif t \nonumber\\
        = & \iint_{Q_T}\!\! \Big\{ \frac{\bq}{\br} \d_t\phi \d_t f^0 + \d_t\psi \cdot \d_t f \Big\}\dif x \dif t + \frac{\bq}{\br}\! \iint_{Q_T}\!\! \d_t\phi \big\{  \frac{\d_t\phi}{2} \dd u - (\d_t\psi\cdot \nabla )\phi \big\}\dif x \dif t. \nonumber\\
        \le& C \Big\{ |u_b| N^2(T) + N^3(T)\Big\}. \label{temp:timeD3}
    \end{align}
    Taking the derivative $\d_t^{k}$ on \eqref{E0} for $k=0,1,$ and multiplying the resultant equality by $\d_t^{k+1}\phi$, we obtain 
    \begin{align}
        |\d_t^{k+1}\phi|^2 = \d_t^{k+1} \phi \big\{ \d_t f^0 - \rho_+ \dd \d_t^k \psi - \d_t^k\big(u\cdot \nabla \phi\big) \big\}. \label{temp:timeD4}
    \end{align}
    Taking the derivative $\d_t^k$ on \eqref{Ei} for $k=0,1$, and taking inner-product on the resultant equality with $\d_t^{k+1} \psi$, we get
    \begin{align}
        &\d_t\Big\{ \frac{\mu}{2\br} |\nabla \d_t^k \psi|^2 + \frac{\mu+\lambda}{2\br}\big( \dd \d_t^k \psi - \frac{\br\bq}{\mu+\lambda} \d_t^k\phi \big)^2 - \frac{\br\bq^2}{2(\mu+\lambda)}|\d_t^k\phi|^2 \Big\}+|\d_t^{k+1}\psi|^2 \nonumber\\
        =& \d_t^{k+1} \psi \cdot \d_t^k f - \bq \d_t^{k+1} \phi \dd \d_t^k \psi \nonumber\\
        &+ \dd\Big\{ \frac{\mu}{\br} \big(\nabla \d_t^k \psi\big)^{\top} \cdot \d_t^{k+1}\psi + \frac{\mu+\lambda}{\br} \d_t^{k+1}\psi \dd \d_t^k \psi - \bq \d_t^k \phi \d_t^{k+1}\psi \Big\}. \label{temp:timeD5}
    \end{align}
    Summing equations \eqref{temp:timeD4}--\eqref{temp:timeD5}, integrating in $(t,x)\in Q_T$, using boundary condition \eqref{dt-bdry} and inequality \eqref{Linf}, we get
    \begin{align*}
        &\int_{\Omega}\!\! \Big\{ \frac{\mu}{2\br} |\nabla \d_t^k \psi|^2 + \frac{\mu+\lambda}{2\br}\big( \dd \d_t^k \psi - \frac{\br\bq}{\mu+\lambda} \d_t^k\phi \big)^2\Big\}\dif x \bigg\vert_{t=0}^{t=T} + \iint_{Q_T}\!\!\! \big| \d_t^{k+1}(\phi,\psi) \big|^2\, \dif x \dif t\\
        =& \frac{\br\bq^2}{2(\mu+\lambda)} \|\d_t^k\phi(t,\cdot)\|^2\big\vert_{t=0}^{t=T} + \iint_{Q_T} \d_t^{k+1}\psi \cdot \d_t^k f\, \dif x \dif t\\
        & +\iint_{Q_T}\!\!\! \d_t^{k+1}\phi \Big\{ \d_t f^0 - (\br+\bq) \dd \d_t^k \psi - \d_t^k\big(u\cdot \nabla \phi\big) \Big\} \dif x \dif t,
    \end{align*}
    which holds true for $k=0,1$. Applying Lemma \ref{lemma:rE}, inequality \eqref{Linf} and the previous estimate \eqref{temp:timeD3} on the above
    equality, we get
    \begin{align*}
        &\int_{\Omega}\!\! \Big\{ \frac{\mu}{2\br} |\nabla \d_t^k \psi|^2 + \frac{\mu+\lambda}{2\br}\big( \dd \d_t^k \psi - \frac{\br\bq}{\mu+\lambda} \d_t^k\phi \big)^2\Big\}\dif x + \frac{1}{2}\iint_{Q_T}\!\!\! \big| \d_t^{k+1}(\phi,\psi) \big|^2\, \dif x \dif t\\
        \le & C \Big\{ N^2(0) + |u_b| N^2(T) + N^3(T) \Big\}.
    \end{align*}
    This completes the proof.
\end{proof}

\subsection{Estimates for angular derivatives}\label{ssec:angular}

\paragraph{} Our aim is to obtain the estimates for spatial derivatives $\|(\nabla^k \psi, \nabla^k \phi)(t)\|$ where $k=1,2,3$. However the boundary condition \eqref{bdry} poses problem for the derivation of parabolic estimates using equations \eqref{Ei}. To resolve this issue, we decompose the gradient operator $\nabla$ into the tangential and normal differential operators, relative to the boundary surface $\d\Omega$. Since the tangential differential operators preserves the boundary condition, it is viable to derive the parabolic estimates for the tangential derivatives of $(\psi,\phi)$. Thus the main purpose of this subsection is to specify the notion of tangential differential operators $\d$ and obtain their corresponding estimates.  

We consider the spherical coordinate system $(r,\theta,\vp)$, where $r\ge 1$ is the radial variable, $\theta\in[0,\pi]$ is the polar angle and $\vp\in[0,2\pi)$ is the azimuthal angle. For the angular derivatives, terms such as $1/\sin\theta$ occur, which form singularity at $\theta=0$ or $\pi$. To circumvent this issue, we employ the following partition of unity
\begin{align*}
    \Omega_V \vcentcolon= \Big\{ x=(x_1,x_2,x_3)^{\top}\in \R^3 \ \Big\vert\ |x|\ge 1, \quad \arccos\big( \frac{x_3}{|x|} \big)\in \big[\frac{\pi}{9},\frac{8\pi}{9}\big]  \Big\},\\
    \Omega_H \vcentcolon= \Big\{ x=(x_1,x_2,x_3)^{\top}\in \R^3 \ \Big\vert\ |x|\ge 1, \quad  \arccos\big( \frac{x_2}{|x|} \big)\in \big[\frac{\pi}{9},\frac{8\pi}{9}\big]  \Big\}.
\end{align*}
It is obvious that $\Omega=\Omega_V\cup \Omega_H$. On each domain $\Omega_V$ or $\Omega_H$, we define the spherical coordinate variables $(r,\theta_V,\vp_V)$ and $(r,\theta_H,\vp_H)$ as 
\begin{alignat*}{4}
    &x_1 = r \cos\vp_V \sin\theta_V, \quad &&x_2 = r\sin\vp_V \sin\theta_V, \quad &&x_3 = r \cos\theta_V, \quad &&\text{for } x\in \Omega_V,\\
    &x_1 = r \cos\vp_H \sin\theta_H, \quad &&x_2 = r \cos\theta_H, \quad &&x_3 = r \sin\vp_H \sin\theta_H, \quad &&\text{for } x\in \Omega_H.
\end{alignat*}
For $x\in \Omega_V$, we also define the spherical unit vectors 
\begin{align*}
    \hr \vcentcolon= \frac{x}{|x|}, \quad \hth_{V} \vcentcolon= \frac{ x_1 x_3 \he_1 +  x_2 x_3 \he_2 - (x_1^2+x_2^2) \he_3 }{|x|\sqrt{x_1^2 + x_2^2}}, \quad \hvp_V \vcentcolon= \frac{-x_2 \he_1 + x_1 \he_2}{\sqrt{x_1^2 + x_2^2}}.
\end{align*}
For $x\in \Omega_H$, unit vectors $(\hr,\hth_H,\hvp_H)$ is defined by exchanging the terms $ x_2 \leftrightarrow x_3$ and $\he_2 \leftrightarrow \he_3$ on the above expression.
Let $\chi_V$, $\chi_H\in \mC^{\infty}(\R^3)$ be non-negative functions satisfying
\begin{subequations}\label{chi}
\begin{gather}
    0\le \chi_V, \, \chi_H \le 1, \quad \chi_{V} + \chi_{V} \ge 1, \quad  \supp\chi_V \subseteq \Omega_V, \quad \supp\chi_H \subseteq \Omega_H,\label{chi1}\\
    |\nabla\chi_H|^2 \le C \chi_H, \qquad |\nabla\chi_V|^2 \le C \chi_V, \qquad |D^k\chi_{H}|+ |D^k\chi_V| \le C,\label{chi2}\\
     \hr\cdot \nabla \chi_H= \hvp_H\cdot \nabla \chi_H = 0, \qquad \hr\cdot \nabla \chi_V= \hvp_V\cdot \nabla \chi_V = 0\label{chi3}
\end{gather}   
\end{subequations}
 for $x\in\R^3$ and $k=1,2,3,4$,
where $C>0$ is a certain generic constant. The detailed construction for $\chi_V$ and $\chi_H$ is explained in Appendix \ref{append:sph}. 
For $x\in\Omega_{V}$ or $\Omega_{H}$,
we define the radial and angular derivatives as follows
\begin{gather*}
    \d_r \vcentcolon= \frac{x}{|x|}\cdot \nabla, \qquad \d_{\theta_V} \vcentcolon= |x| \hth_{V} \cdot \nabla, \qquad \d_{\vp_V} \vcentcolon= \sqrt{x_1^2 + x_2^2} \,\, \hvp_V \cdot \nabla,\\
    \d_{\theta_H} \vcentcolon= |x| \hth_{H} \cdot \nabla, \qquad \d_{\vp_H} \vcentcolon= \sqrt{x_1^2 + x_3^2} \,\, \hvp_H \cdot \nabla.
\end{gather*}
For simplicity, we use the following notations
\begin{equation*}\label{partial}
     \d \in \{\d_{\theta_V}, \d_{\vp_V}, \d_{\theta_H}, \d_{\vp_H}\},  \qquad D\in \{\d_r, \d_{\theta_V}, \d_{\varphi_V}, \d_{\theta_H}, \d_{\vp_H}\}.
\end{equation*}
In what follows, we only derive estimates for one of $\Omega_V$ and $\Omega_H$ since the other follows exactly from the same argument. To simplify, we use the notations
\begin{subequations}\label{chi-abb}
    \begin{gather}
    \chi \vcentcolon = \chi_{V} \ \text{or} \ \chi_{H}, \qquad \tilde{\Omega} \vcentcolon = \Omega_V \ \text{or} \ \Omega_H, \\ 
    (\theta,\vp) \vcentcolon = (\theta_V,\vp_V) \ \text{or} \ (\theta_H,\vp_H), \qquad \d \in \big\{ \d_\theta, \d_\vp \big\}, \qquad D\in \big\{ \d_r, \d_\theta, \d_\vp \big\},\\
    \tilde{Q}_T \vcentcolon = [0,T]\times\Omega_V \ \text{or} \ [0,T]\times\Omega_H, \qquad \tilde{\Gamma}_T \vcentcolon = [0,T]\times\d\Omega_V \ \text{or} \ [0,T]\times\d\Omega_H.
\end{gather}
\end{subequations}
We remark that $\d$ preserves boundary condition \eqref{bdry} in the sense that
\begin{equation}\label{dpsi-bdry}
    \chi \d^k\psi \vert_{x\in\d\Omega} =0 \qquad \text{for } \ \d \in\{ \d_{\theta}, \d_{\vp} \} \ \text{ and } \ k=1,2,3.
\end{equation}
\begin{lemma}\label{lemma:dpsi}
    Suppose $(\phi,\psi)\in X(T_0)$ is a solution to the system \eqref{LE}, where $T_0>0$ is the time given in \eqref{continuity}--\eqref{priorirho}. Then there exists a constant $C=C(\br,\gamma,K,\mu,\lambda)>0$ such that for arbitrary $\delta\in(0,1)$ and $T\in[0,T_0]$,
    it holds that
    \begin{align*}
        &\sup\limits_{0\le t\le T} \big\|\sqrt{\chi}(\d^k \psi,\d^k\phi)(t,\cdot)\big\|^2  + \int_{0}^T\!\!\! \| \sqrt{\chi} D \d^{k}\psi (t,\cdot) \|^2\, \dif t\\ \le& \frac{C}{\delta^7} \Big\{ N^2(0) + |u_b| N^2(T) +N^3(T) \Big\} + C \delta N^2(T),
    \end{align*}
    where $k=1,2,3$.
\end{lemma}
\begin{proof}
    Taking the derivative $\d^k$ on equation \eqref{E0}
    and multiplying both sides of the resultant equality by $\chi \d^k \phi$, we get
    \begin{align*}
        \frac{\bq}{2\br} \d_t (\chi |\d^k \phi|^2 ) + \bq \chi \d^k \phi\, \dd \d^k \psi = \bq \chi \d^k \phi [\dd,\d^k] \psi + \frac{\bq}{\br} \chi \d^k\phi \d^k \big\{ f^0 - (u\cdot \nabla) \phi \big\}.
    \end{align*}
    Taking the derivative $\d^k$ on equation \eqref{Ei}
    and multiplying both sides of the resultant equality by $\chi \d^k\psi$, we get
    \begin{align*}
        &\frac{1}{2}\d_t(\chi |\d^k \psi|^2) - \frac{\mu}{\br} \chi \d^k \psi \cdot \Delta \d^k \psi -\frac{\mu+\lambda}{\br} \chi  (\d^k\psi \cdot \nabla) \dd \d^k \psi + \bq \chi  (\d^k\psi \cdot \nabla) \d^k \phi \\
        =& \bq \chi \d^k \psi \cdot [\nabla, \d^k]\phi + \frac{\mu}{\br} \chi \d^k\psi \cdot [\d^k,\Delta] \psi + \frac{\mu+\lambda}{\br} \chi \d^k\psi \cdot \big[ \d^k, \nabla \dd \big] \psi + \chi \d^k \psi \cdot \d^k f.
    \end{align*}
    Summing 
    the above two equations, integrating
    the resultant
    equality
    in $(t,x)\in Q_T \vcentcolon = [0,T]\times \Omega$ with integrating by parts 
    using
    boundary condtion
    \eqref{dpsi-bdry}, we have that
    \begin{align}
        &\int_{\Omega} \frac{\chi}{2} \big\{  |\d^k \psi|^2 + \frac{\bq}{\br} |\d^k \phi|^2 \big\}(t,x) \, \dif x\Big\vert_{t=0}^{t=T}\nonumber\\
        &+ \iint_{Q_T} \!\frac{\chi}{\br} \Big\{ \mu |\nabla \d^k \psi|^2  + (\mu+\lambda)|\dd \d^k \psi|^2 \Big\} \, \dif x \dif t\nonumber\\
        =& \iint_{Q_T}\!\! \Big\{ \bq \d^k\phi (\d^k\psi \cdot \nabla) \chi - \frac{\mu}{\br} \d^k \psi^{\top} \cdot \nabla \d^k \psi  \cdot \nabla \chi - \frac{\mu+\lambda}{\br} \dd \d^k \psi ( \d^k \psi \cdot\nabla) \chi   \Big\} \, \dif x \dif t \nonumber\\
        & + \iint_{Q_T} \chi \Big\{ \bq \d^k \phi \cdot [\dd , \d^k]\psi +  \bq \d^k \psi \cdot \big[ \nabla, \d^k \big]\phi \Big\} \, \dif x \dif t\nonumber\\
        &+ \iint_{Q_T} \chi \Big\{ \d^k \psi \cdot \d^k f + \frac{\bq}{\br} \d^k \phi \d^k(f^0 - u\cdot \nabla\phi) \Big\} \, \dif x \dif t  \nonumber\\
        &+ \frac{1}{\br}\iint_{Q_T} \chi \d^k \psi \cdot \Big\{ \mu [\d^k,\Delta] \psi + (\mu+\lambda) \big[ \d^k, \nabla \dd \big] \psi \Big\} \, \dif x \dif t =\vcentcolon \sum_{i=1}^4 \mI_i.\label{temp:dpsi1}
    \end{align}
    By the property of cut-off function $\chi$ in \eqref{chi}, we have for an arbitrary $\epsilon\in(0,1)$
    \begin{align*}
        |\mI_1 | 
        \le& \iint_{Q_T}\!\! \frac{\chi}{4\br} \Big\{ \mu |\nabla \d^k \psi|^2 + (\mu+\lambda) |\dd \d^k \psi|^2 \Big\}\, \dif x \dif t  + \frac{\epsilon}{3} N^2(T) + \frac{C}{\epsilon} \iint_{\tilde{Q}_T}\!\! |\d^k\psi|^2 \, \dif x \dif t.
    \end{align*}
    Using Proposition \ref{prop:dcom}, we also get
    \begin{align*}
        |\mI_2| 
        \le & \frac{\epsilon}{2} N^2(T) + \frac{C}{\epsilon} \sum_{m=0}^{k-1} \iint_{Q_T} \chi |\nabla \d^{m}\psi|^2 \, \dif x \dif t.
    \end{align*}
    For the term $\mI_3$, we recall the definition \eqref{fi}. Then integrating by parts, applying Lemma \ref{lem:st} and inequality \eqref{Linf}, we get for $k=1,2,3$,
    \begin{align*}
        &\bigg|\iint_{Q_T} \chi \d^k \psi \cdot \d^k (f-\tf) \, \bigg| = \bigg|\sum_{i=1}^3 \iint_{Q_T} \chi \d^k \psi^i \d^k \Big( \frac{\phi+\tr-\br}{\br(\phi+\tr)} (L\psi)^i \Big)\bigg| \\
        =& \bigg|\sum_{i=1}^3\! \iint_{Q_T}\!\!\! \big\{\d \chi \d^k \psi^i + \chi \d^{k+1} \psi^i  \big\} \d^{k-1}\! \Big( \frac{\phi+\tr-\br}{\br(\phi+\tr)} (L\psi)^i \Big) \bigg| \!\le\! C \big\{ |u_b| N^2(T) + N^3(T) \big\}.
    \end{align*}
    Thus by the above estimate, Lemma \ref{lem:st}, and inequality \eqref{Linf}, it follows that
    \begin{align*}
        |\mI_3| = & \bigg| \iint_{Q_T} \chi \d^k \psi \cdot \big\{ \d^k f + \d^k(f-\tf) \big\}  + \chi \frac{\bq}{\br} \d^k \phi \d^k(f^0-u\cdot \nabla \phi)  \dif x \dif t \bigg| \\ 
        \le & C N^3(T) + C|u_b| N^2(T). 
    \end{align*}
    Applying Propositions \ref{prop:hardy} and \ref{prop:dcom} we obtain that
    \begin{align*}
        |\mI_4| 
        \le& C \iint_{Q_T} \chi |\d^k\psi| \Big\{  \sum_{m=1}^{k+1}|\d^m \psi| + \frac{|\psi|^2}{|x|^4} + \sum_{m=0}^{k-1} |D^2 \d^m \psi| \Big\} \, \dif x \dif t\\ 
        \le& \frac{\epsilon}{3} N^2(T) + \frac{C}{\epsilon}\iint_{Q_T} \chi |\d^k \psi|^2 \, \dif x \dif t.
    \end{align*}
    Substituting
    the estimates
    for
    $\mI_1$--$\mI_4$ in
\eqref{temp:dpsi1}, we have that
    \begin{align}
        &\int_{\Omega}\! \chi\big\{  |\d^k \psi|^2 + \frac{\bq}{\br} |\d^k \phi|^2 \big\}(T,x) \, \dif x + \iint_{Q_T} \!\chi \Big\{ \mu |\nabla \d^k \psi|^2  + (\mu+\lambda)|\dd \d^k \psi|^2 \Big\} \, \dif x \dif t\nonumber\\
        \le & C \big\{ N^2(0) + |u_b|N^2(T) + N^3(T) \big\} + \epsilon N^2(T) + \frac{C}{\epsilon} \sum_{m=0}^{k-1} \iint_{\tilde{Q}_T} |\nabla \d^m \psi|^2 \, \dif x \dif t,\label{temp:dpsi2}
    \end{align}
    where $k=1,2,3$ and $\epsilon\in(0,1)$. For an arbitrary $\delta\in(0,1)$, let $\epsilon= \delta^4$ and $k=1$ in \eqref{temp:dpsi2}. Then 
    we have from Lemma \ref{lemma:rE} that
    \begin{align}
        &\int_{\Omega}\! \chi \big\{  |\d \psi|^2 + |\d \phi|^2 \big\}(T,x) \, \dif x + \iint_{Q_T}\! \chi \Big\{ \mu |\nabla \d \psi|^2  + (\mu+\lambda)|\dd \d \psi|^2 \Big\} \, \dif x \dif t\nonumber\\
        \le & \frac{C}{\delta^4}\Big\{ N^2(0)+ |u_b| N^2(T) + N^3(T) \Big\} + \delta^4 N^2(T).\label{temp:dpsi3}
    \end{align}
    Next we let
    $\epsilon = \delta^2$ and $k=2$ in \eqref{temp:dpsi2}. Applying Lemma \ref{lemma:rE} and substituting \eqref{temp:dpsi3} into the resultant inequality, we get
    \begin{align}
        &\int_{\Omega}\! \chi \big\{  |\d^2 \psi|^2 + |\d^2 \phi|^2 \big\}(T,x) \, \dif x + \iint_{Q_T} \!\! \chi\Big\{ \mu |\nabla \d^2 \psi|^2  + (\mu+\lambda)|\dd \d^2 \psi|^2 \Big\} \, \dif x \dif t\nonumber\\
        \le& \frac{C}{\delta^6} \Big\{ N^2(0) + |u_b| N^2(T) +  N^3(T) \Big\} + C \delta^2 N^2(T). \label{temp:dpsi4}
    \end{align}
    Finally, let $\epsilon = \delta$ and $k=3$ in \eqref{temp:dpsi2}. Applying Lemma \ref{lemma:rE} and substituting \eqref{temp:dpsi3}--\eqref{temp:dpsi4} into the resultant inequality, we get
    \begin{align*}
        &\int_{\Omega} \! \chi \big\{  |\d^3 \psi|^2 + |\d^3 \phi|^2 \big\}(T,x) \, \dif x + \iint_{Q_T} \!\! \chi\Big\{ \mu |\nabla \d^3 \psi|^2  + (\mu+\lambda)|\dd \d^3 \psi|^2 \Big\} \, \dif x \dif t\\ 
        \le & \frac{C}{\delta^7} \Big\{ N^2(0) + |u_b| N^2(T) +N^3(T) \Big\} + C \delta N^2(T) .
    \end{align*}
    This concludes the proof.
\end{proof}

\begin{corollary}\label{corol:dkdtphi}
    Suppose the same assumption in Lemma \ref{lemma:dpsi} holds. Then there exists a constant $C=C(\br,\gamma,K,\mu,\lambda)>0$ such that for arbitrary $\delta\in(0,1)$ and $T\in[0,T_0]$, it holds that
    \begin{align*}
        &\iint_{Q_T} |D_t \phi|^2 \, \dif x \dif t \le C\Big\{ N^2(0) + N^3(T) + |u_b| N^2(T) \Big\},\\
        &\iint_{Q_T} \chi |\d^{k}D_t \phi|^2 \, \dif x \dif t \le \frac{C}{\delta^7}\Big\{ N^2(0) + |u_b| N^2(T) + N^3(T)  \Big\} + C \delta N^2(T),
    \end{align*}
    where $k=0,1,2,3$.
\end{corollary}
\begin{proof}
    The equation \eqref{ns1} yields $ D_t \phi = - \br \dd \psi + f^0$. Integrating this equation in $(t,x)\in Q_T$, applying Lemmas \ref{lemma:rE} and inequality \eqref{Linf}, we obtain the first estimate. Next, taking the derivative $\d^k$ on \eqref{E0} yields $\d^k D_t \phi = - \br \dd \d^k \psi - \br [\d^k,\dd]\psi + \d^k f^0$. Integrating this equation in $(t,x)\in Q_T$, applying Lemmas \ref{lemma:rE}, \ref{lemma:dpsi}, Proposition \ref{prop:dcom}, and inequality \eqref{Linf}, we obtain the second estimate.
\end{proof}

\subsection{Damping estimates for \texorpdfstring{$\d_r\phi$}{Drphi}}\label{ssec:damp}
\paragraph{} In this subsection, we derive estimates involving the radial derivative of $\phi$. This is obtained by exploiting the coupling structure of the continuity equation \eqref{E0} and momentum equation \eqref{Ei}. More specifically, we show that $\d_r\phi$ solves a first-order hyperbolic equation via the following procedure: first, we take $\d_r$ on \eqref{E0} to get
\begin{align}\label{drdivpsi}
    -\br \d_r \dd \psi =& \d_r D_t \phi - \d_r f^0,
\end{align}
where $D_t\vcentcolon = \d_t + (u\cdot\nabla)$. Next, taking the radial component of \eqref{Ei}, we obtain 
\begin{equation*}
    - \frac{2\mu+\lambda}{\br} \d_r\dd \psi + \bq \d_r\phi = \hr \cdot (f-\d_t\psi) + \frac{\mu}{\br}\big\{ \hr \cdot \Delta \psi - \d_r\dd \psi \big\}.
\end{equation*}
Substituting \eqref{drdivpsi} in the above equation yields
\begin{equation}
     \d_r D_t \phi + \beta \d_r \phi = g, \label{drphi}
\end{equation}
where 
\begin{equation*}
    \beta\vcentcolon= \frac{\br^2 \bq}{2\mu+\lambda}>0, \qquad g \vcentcolon= \d_r f^0 + \frac{\beta}{\bq} \hr\cdot (f-\d_t\psi) + \frac{\mu\beta}{\br\bq}\big\{ \hr\cdot \Delta\psi -\d_r \dd \psi \big\}.
\end{equation*}
Furthermore, computing the commutator $[\d_r, D_t]$, we also obtain that
\begin{equation}\label{drphi-2}
    D_t \d_r \phi + \beta \d_r\phi = g + \d_r u \cdot \nabla \phi .
\end{equation}
The term $\beta \d_r\phi$ on the left hand side indicates the dissipative damping effect on $\d_r\phi$.
\begin{lemma}\label{lemma:damping}
    Suppose $(\phi,\psi)\in X(T_0)$ is a solution to the system \eqref{LE}, where $T_0>0$ is the time given in \eqref{continuity}--\eqref{priorirho}. Then there exists a constant $C=C(\br,\gamma,K,\mu,\lambda)>0$ such that for an arbitrary $T\in[0,T_0]$, it holds that 
    \begin{align*}
        &\sup\limits_{0\le t\le T}\big\|\sqrt{\chi}\d_r^{l+1}\d^k \phi \big\|^2 + \int_{0}^T\!\!\big\|\sqrt{\chi}\d_r^{l+1}\d^k (\phi, D_t\phi) \big\|^2 \, \dif t + |u_b|  \iint_{\tilde{\Gamma}_T} \!\! |\d_r^{l+1}\d^k \phi|^2 \, \dif x \dif t \\
        \le& C \big\{ N^2(0) + |u_b| N^2(T) + N^3(T) \big\}  + C \!\! \int_{0}^T\!\!\! \Big\{ \sum_{j=0}^{k+1} \big\| \sqrt{\chi} D \d^{j} \psi \big\|_{l}^2 +  \big\|\sqrt{\chi}\d_r^l \d^k \d_t\psi\big\|^2 \Big\} \,\dif t,
    \end{align*}
    where $k$, $l$ are non-negative integers such that $k+l=0,1,2$.
\end{lemma}
\begin{proof}
    Taking the derivative $\d_r^l \d^k$ on \eqref{drphi} for $k+l=0,1,2$, multiplying both sides by $\chi \rho \d_r^{l+1}\d^k D_t \phi$, and using continuity equation \eqref{ns1}, we get
    \begin{align*}
        &\d_t \big( \frac{\beta}{2} \chi \rho |\d_r^{l+1}\d^k \phi|^2 \big) + \dd \big( \frac{\beta}{2} \chi \rho u |\d_r^{l+1}\d^k \phi|^2 \big) + |\d_r^{l+1}\d^k D_t \phi|^2 \\
        =& \beta \chi \rho \d_r^{l+1}\d^k \phi [ D_t, \d_r^{l+1} \d^k ] \phi +\chi \rho \d_r^l\d^k g \d_{r}^{l+1} \d^k D_t\phi + \beta\frac{\rho}{2} |\d_r^{l+1}\d^k\phi|^2 (u\cdot \nabla) \chi.  
    \end{align*}
    In addition, taking the derivative $\d_r^l \d^k$ on \eqref{drphi-2}, multiplying both sides by $ \chi \rho \d_r^{l+1} \d^k \phi$, and using continuity equation \eqref{ns1}, we get
    \begin{align*}
        &\d_t \big( \frac{1}{2} \chi \rho |\d_r^{l+1} \d^k \phi|^2 \big) + \dd \big( \frac{1}{2} \chi \rho u |\d_r^{l+1} \d^k \phi|^2 \big)  + \beta \chi \rho \big|\d_r^{l+1}\d^k \phi\big|^2 \\ 
        =& \chi \rho \d_r^{l+1}\d^k \phi \Big\{ \d_r^l \d^k g + \d_r^l \d^k \big( \d_r u \cdot \nabla \phi \big) + [ D_t \d_r \, , \d_r^{l} \d^k ] \phi \Big\} + \frac{\rho}{2} |\d_r^{l+1}\d^k \phi|^2  (u\cdot \nabla) \chi . 
    \end{align*}
    Summing the above equations, integrating the resultant equality in $(t,x)\in[0,T]\times \Omega$, it follows from boundary condition \eqref{bdry} and Cauchy-Schwarz's inequality that
    \begin{align}
        &\frac{\beta+1}{2} \int_{\Omega} \chi \rho |\d_r^{l+1} \d^k \phi|^2(t,x)\, \dif x\Big\vert_{t=0}^{t=T} + \frac{\beta+1}{2}\int_{0}^T\!\!\!\! \int_{\d\Omega}\!\!\! |u_b| \chi \rho |\d_r^{l+1}\d^k \phi|^2 (t,y) \, \dif y \dif t  \nonumber\\
        &+ \frac{1}{2} \iint_{Q_T} \chi \rho \Big\{ |\d_r^{l+1}\d^k D_t \phi|^2 + \beta|\d_r^{l+1} \d^k \phi|^2 \Big\}\, \dif x \dif t \nonumber\\
        \le& C\iint_{Q_T} \chi \Big\{ \big| [D_t\,, \d_r^{l+1} \d^k] \phi \big|^2 + \big| [ D_t \d_r \, , \d_r^{l} \d^k ] \phi\big|^2  \Big\}\, \dif x \dif t \nonumber\\
        &+ C\iint_{Q_T} \chi \Big\{ \big|\d_r^l \d^k ( \d_r\psi \cdot \nabla \phi ) \big|^2 + \big|\d_r^l \d^k ( \d_r\tu \cdot \nabla \phi ) \big|^2  \Big\}\, \dif x \dif t \nonumber\\
        & + C\!\!\iint_{Q_T}\!\! \chi \big\{ |\d_r^{l+1}\d^k f^0|^2 + |f|^2  \big\} \, \dif x \dif t + C\!\!\iint_{Q_T}\!\! \chi \big| \d_r^l \d^k (\hr \cdot \Delta \psi - \d_r \dd \psi ) \big|^2 \, \dif x \dif t \nonumber\\
        & +C \iint_{Q_T} |\d_r^{l+1}\d^k \phi|^2  |(u\cdot \nabla) \chi| \, \dif x \dif t + C\iint_{Q_T}\!\! \chi |\d_r^l \d^k\d_t\psi|^2 \, \dif x \dif t = \vcentcolon \sum_{i=1}^6 \mJ_i.\label{temp:damping1}
    \end{align}
    Using the inequality \eqref{Linf}, Proposition \ref{prop:dcom}, and the fact that $\d_t$, $\d_r$, $\d$ commute with each other, we obtain
    \begin{align*}
        |\mJ_1| =& \bigg| \iint_{Q_T}\!\!  \chi \Big\{ \big([ (u\cdot \nabla) \, , \d_r^{l+1} \d^k ] \phi\big)^2 + \big([ (u\cdot \nabla) \d_r \, , \d_r^{l} \d^k ] \phi\big)^2  \Big\}\, \dif x \dif t \bigg|\\ 
        \le&  C |u_b|^2 N^2(T) + C N^4(T),\\
        |\mJ_2| \le&  C |u_b|^2 N^2(T) + N^4(T).
    \end{align*}
    By the definition \eqref{f0}--\eqref{fi} and inequality \eqref{Linf}, we obtain 
    \begin{align*}
        |\mJ_3| = C\bigg|\iint_{Q_T} \chi \big\{ |\d_r^{l+1}\d^k f^0|^2 + |f|^2  \big\} \, \dif x \dif t \bigg| \le C |u_b| N^2(T) + C N^3(T).
    \end{align*}
    Owing to Proposition \ref{prop:rrcancel}, it follows that
    \begin{align*}
        |\mJ_4| =& C \bigg|\iint_{Q_T} \chi |\d_r^{l}\d^k(\hr\cdot\Delta\psi -\d_r \dd \psi )|^2 \, \dif x \dif t \bigg| \le C \int_{0}^T \sum_{j=0}^{k+1} \big\|\sqrt{\chi} D \d^{j} \psi \big\|_{l}^2\,\dif t.
    \end{align*}
    Using Lemma \ref{lem:st} and inequality \eqref{Linf}, we get
    \begin{align*}
        |\mJ_5| = C \bigg| \iint_{Q_T} |\d_r^{l+1}\d^k \phi|^2  |(\psi+\tu)\cdot \nabla \chi| \, \dif x \bigg| \le C |u_b| N^2(T) + C N^3(T). 
    \end{align*}
    Substituting the estimates 
    for $\mJ_1$--$\mJ_5$ in
    \eqref{temp:damping1}, we obtain the desired estimate.
\end{proof}

\subsection{Estimates from elliptic and Stokes systems}\label{ssec:stokes}
\paragraph{} Using the previous results derived in Sections \ref{ssec:ETD}--\ref{ssec:damp}, we obtain the estimates containing full spatial derivatives $\|(\nabla^k\psi,\nabla^k\phi)\|$ for $k=1,2,3$. This is achieved by reformulating equations \eqref{E0}--\eqref{Ei} and applying elliptic and Stokes-type estimates. To illustrate this, we first present the two known results for solutions to the elliptic system and Stokes' equations.   

Let $\nu_1,\,\nu_2 >0$. Given a vector-valued function $\mathcal{R}\vcentcolon\Omega \to \R^3$, which belongs to $H^m(\Omega)$ for $m=0,1$. Suppose $\zeta\vcentcolon \Omega\to\R^3$ is a solution to the elliptic system:
\begin{subequations}\label{elliptic}
    \begin{align}
    - \nu_1 \Delta \zeta^i - \nu_2 \d_{x_i} \dd \zeta =& \mathcal{R}^i  &&\text{for } \ x\in\Omega \ \text{ and } \  i=1,2,3, \\
    \zeta(x) =& 0  && \text{for } \ x\in\d\Omega.
\end{align}
\end{subequations}

\begin{lemma}\label{lemma:elliptic}
    There exists a constant $C=C(\nu_1,\nu_2)>0$ such that for $k=2,3$,
    \begin{align*}
        \|D^k \zeta\| \le C \big\{ \|\mathcal{R}\|_{k-2} + \|\zeta\| \big\}.
    \end{align*}
\end{lemma}
The above estimate is well-known and its derivation is presented in \cite{ADN1964}. The $L^2$ norm on $\zeta$ contained in the right hand side is due to the unboundedness of the exterior domain $\Omega$.  

In addition, suppose $\alpha\vcentcolon\d\Omega\to \R^3$, $h^0 \vcentcolon \Omega \to \R$, and $h\vcentcolon \Omega \to \R^3$ are functions which belong to the suitable Sobolev spaces. Let $w\vcentcolon \Omega \to \R$ and $v\vcentcolon \Omega \to \R^3$ be solution to the Stokes-type equations:
\begin{subequations}\label{stokes}
    \begin{align}
        \br \dd v  &= h^0  &&\text{for } \ x\in \Omega,\label{St0}\\
        -\mu \Delta v^i + \br \bq \d_{x_i}w &= h^i  &&\text{for } \ x\in \Omega \ \text{ and } \ i=1,2,3, \label{Sti}\\
        v(x)&= \alpha(x) &&\text{for } \ x\in\d\Omega. \label{St-bdry}
    \end{align}
\end{subequations}
\begin{lemma}\label{lemma:stok}
    There exists $C=C(\br,\gamma,K,\mu,\lambda)>0$ such that for $k=2,3,4$,
    \begin{align*}
        \|D^k v\|^2 + \|D^{k-1} w\|^2 \le C \big\{ \|h^0\|_{k-1}^2 + \|h\|_{k-2}^2 + \|\alpha\|_{H^{k-1/2}(\d\Omega)}^2 + \|Dv\|^2 \big\}.
    \end{align*}
\end{lemma}
The proof for the above lemma is presented in \cite{m-n83} (See Lemma 4.3). We also note that for the incompressible case: $h^0\vcentcolon = 0$, the system \eqref{stokes} reduces to the well-known Stokes equations, and the similar estimates are derived in \cite{fin1965,hey1979}.  
\begin{lemma}\label{lemma:stokes}
    There exists a constant $C=C(\br,\gamma,K,\mu,\lambda)>0$ such that for $l=0,1,2$, 
    \begin{align*}
        \|D^{2+l} \psi\|^2 + \| D^{l+1} \phi \|^2 \le C \Big\{ \|\d_t \psi\|_{l}^2 + \| D_t \phi \|_{l+1}^2 + \|D \psi\|^2 +\|f^0\|_{l+1}^2 + \|f\|_{l}^2 \Big\}.
    \end{align*}
    Moreover, for $k\ge 1$ and $k+l=1,2$,
    \begin{align*}
        \big\| \chi D^{2+l} \d^k \psi \big\|^2 + \big\| \chi D^{1+l} \d^k \phi \big\|^2
        \!\le & C \Big\{ \big\| \chi \d^k D_t \phi \big\|_{l+1}^2 \!\!+ \big\| \chi \d_t \d^k \psi \big\|_{l}^2 \! + \big\| f^0 \big\|_{k+l+1}^2\!\! + \big\| f \big\|_{k+l}^2  \Big\}\\
        &+ C \Big\{ \sum_{m=1}^{k}\big\| \d^m \psi \big\|_{H^{l+1}(\tilde{\Omega})}^2 + \sum_{m=0}^{k-1} \big\| D \d^m \phi \big\|_{H^{l}(\tilde{\Omega})}^2 \Big\},
    \end{align*}
    where $\tilde{\Omega}$ is the domain defined in \eqref{chi-abb}.
\end{lemma}
\begin{proof}
    Applying the operator $\chi\d^k$ on the equations \eqref{E0}--\eqref{Ei}, we have that $(v,w)= (\chi\d^k\phi,\chi\d^k\psi)$ solves the Stokes-type system \eqref{stokes} with $\chi\d^k\psi \vert_{\d\Omega} = 0 = \alpha$, and
    \begin{align*}
        h^0 \vcentcolon=& \chi \d^k f^0 - \chi \d^k D_t \phi + \br \big(\d^k \psi \cdot \nabla \big)\chi,\\
        h^{i} \vcentcolon=& \br \chi \d^k f^i - \br \chi \d_t \d^k \psi^i + \frac{\mu+\lambda}{\br} \chi \d^k \d_{x_i}(f^0-D_t\phi)\\
        &+ \mu \chi [\d^k,\Delta]\psi^i - \br\bq \chi [\d^k, \d_{x_i}] \phi + \mu \Delta\chi \d^k \psi^i + 2\mu (\nabla \chi \cdot \nabla) \d^k \psi^i,
    \end{align*}
    By \eqref{chi}, we have $|(\nabla \chi \cdot \nabla) F| \le C \sqrt{\chi} | \d_{\theta} F |$. Then by Lemma \ref{lemma:stokes} and Propositions \ref{prop:dcom}, we have
    \begin{align*}
        &\big\| \chi D^{2+l} \d^k \psi \big\|^2 + \big\| \chi D^{1+l} \d^k \phi \big\|^2\\ 
        \le & C \Big\{ \big\| \chi \d^k D_t \phi \big\|_{l+1}^2 + \big\| \chi \d_t \d^k \psi \big\|_{l}^2 + \big\| f^0 \big\|_{k+l+1}^2 + \big\| f \big\|_{k+l}^2  \Big\}\\
        &+ C \Big\{ \sum_{m=1}^{k}\big\| \d^m \psi \big\|_{H^{l+1}(\tilde{\Omega})}^2 + \sum_{m=0}^{k-1} \big\| D \d^m \phi \big\|_{H^{l}(\tilde{\Omega})}^2 \Big\}.
    \end{align*}
    This proves the second inequality. The first statement is shown by the exact same argument except we set $v=\psi$, $w=\phi$ in Lemma \ref{lemma:stok}.
\end{proof}
\begin{lemma}\label{lemma:Htwo}
    There exists a constant $C=C(\br,\gamma,K,\mu,\lambda)>0$ such that for $\delta\in(0,1)$ and $T\in[0,T_0]$,
    \begin{align*}
        &\sup\limits_{0\le t\le T}\|(\phi, \psi)\|_2^2 + \int_{0}^T\!\! \big\{ \| D\phi \|_1^2 + \|( D \psi, \d_t \phi) \|_2^2  \big\}\dif t \\
        \le& \frac{C}{\delta^7}\Big\{ N^2(0) + |u_b|N^2(T) + N^3(T) \Big\} + C \delta N^2(T). 
    \end{align*}
\end{lemma}
\begin{proof}
    Throughout this proof, we use the abbreviation
    \begin{equation}\label{Mdelta}
        M_{\delta}(T) \vcentcolon= \frac{C}{\delta^7}\Big\{ N^2(0) + |u_b|N^2(T) + N^3(T) \Big\} + C \delta N^2(T),
    \end{equation}
    where $C=C(\br,\gamma,K,\mu,\lambda)>0$ is a constant which is independent of $T\in[0,T_0]$, $|u_b|$, $\delta\in(0,1)$ and the initial data $(\rho_0,u_0)$.
    
    First, set $(k,l)=(0,0)$ and $(k,l)=(1,0)$ in Lemma \ref{lemma:damping}. Then using Lemmas \ref{lemma:rE}, \ref{lemma:timeD}, \ref{lemma:dpsi}, and Corollary \ref{corol:dkdtphi}, we get
    \begin{align}\label{temp:Htwo1}
        &\sup\limits_{0\le t\le T}\! \big\|\big( D \phi, \sqrt{\chi} D\d\phi\big)\big\|^2 \! + \!\! \int_{0}^T\!\!\! \big\| \big(\d_r\phi, \sqrt{\chi}\d_r \d \phi,D D_t \phi, \sqrt{\chi} D \d D_t \phi\big) \big\|^2 \dif t  
        \le M_\delta(T).
    \end{align}
    Let $(k,l)=(0,0)$ in Lemma \ref{lemma:stokes}. Then by Lemma \ref{lemma:timeD}, Corollary \ref{corol:dkdtphi}, and \eqref{temp:Htwo1},
    \begin{align}
        &\int_{0}^T \!\! \big\{ \big\|D^2\psi\big\|^2 + \|D\phi\|^2 \big\}\dif t 
        \le M_{\delta}(T).\label{temp:Htwo2}
    \end{align}
    Let $(k,l)=(1,0)$ in Lemma \ref{lemma:stokes}. Then by Lemmas \ref{lemma:rE}, \ref{lemma:timeD}, \ref{lemma:dpsi}, and \eqref{temp:Htwo1}--\eqref{temp:Htwo2},
    \begin{align}
        &\int_{0}^T\!\! \big\{ \big\| \chi D^2 \d \psi \big\|^2 + \big\| \chi D\d \phi \big\|^2 \big\}\dif t 
        \le M_\delta(T).\label{temp:Htwo3}
    \end{align}
    Let $(k,l)=(0,1)$ in Lemma \ref{lemma:damping}. Then using Lemmas \ref{lemma:rE}, \ref{lemma:timeD}, \ref{lemma:dpsi}, Corollary \ref{corol:dkdtphi}, and \eqref{temp:Htwo2}--\eqref{temp:Htwo3}, we obtain
    \begin{align}
        &\sup\limits_{0\le t \le T} \| D^2 \phi \|^2 + \int_{0}^T\!\! \Big\{ \| \d_r^2 \phi \|^2 + \| D^2 D_t \phi \|^2 \Big\}\dif t 
        \le M_{\delta}(T).\label{temp:Htwo4}
    \end{align}
    Let $(k,l)=(0,1)$ in Lemma \ref{lemma:stokes}. Then using Lemmas \ref{lemma:rE}, \ref{lemma:timeD}, Corollary \ref{corol:dkdtphi}, inequalities \eqref{temp:Htwo1} and \eqref{temp:Htwo4}, we have
    \begin{align*}
        &\int_{0}^T \big\{ \big\| D^3 \psi \big\|^2 + \big\| D^2 \phi \big\|^2  \big\}\dif t 
        \le  M_{\delta}(T).
    \end{align*}
    Next, according to equation \eqref{Ei}, we set $\zeta = \psi$ in Lemma \ref{lemma:elliptic}. Then by Lemma \ref{lemma:timeD}, and \eqref{temp:Htwo1}, it follows that
    \begin{align*}
        \sup\limits_{0\le t\le T}\| D^2 \psi \|^2 \le C \sup\limits_{0\le t\le T}\Big\{ \|\d_t \psi\|^2 + \|D\phi\|^2 + \|f\|^2 \Big\} \le  M_{\delta}(T).
    \end{align*}
    Since $D^2 \d_t \phi = D^2 D_t\phi - D^2 u \cdot \nabla \phi - 2 Du\cdot \nabla \phi - u\cdot D^2 \nabla \phi $, \eqref{Linf} and \eqref{temp:Htwo4} implies that $$\int_{0}^T \| \d_t\phi \|_2^2 \, \dif t \le M_{\delta}(T).$$
    This completes the proof.
\end{proof}

\begin{lemma}\label{lemma:Hthree}
    There exists a constant $C\!=\!C(\br,\gamma,K,\mu,\lambda)>\!0$ such that for $\delta\in(0,1)$ and $T\in[0,T_0]$,
    \begin{align*}
        &\sup\limits_{0\le t\le T}\big\|(D^3\phi, D^3\psi)\big\|^2 + \int_{0}^T\!\! \Big\{ \big\|\big( D^3\phi, D^3 \d_t \phi\big) \big\|^2 + \|D^4 \psi\|^2 \Big\}\dif t\\
        \le& \frac{C}{\delta^7}\Big\{ N^2(0) + |u_b| N^2(T) + N^3(T) \Big\} + C \delta N^2(T). 
    \end{align*}
\end{lemma}
\begin{proof}
    For simplicity, we use the same notation $M_{\delta}(T)$ as defined in \eqref{Mdelta}. Taking the temporal derivative $\d_t$ on the equation \eqref{Ei}, we verify that $\zeta = \d_t \psi$ solves the elliptic system \eqref{elliptic}, with $\nu_1 = \mu$, $\nu_2 = \mu+\lambda$ and
    \begin{equation*}
        \mathcal{R}^i = \br \d_t f^i - \br \d_t^2 \psi^i - \br(\d_t\psi \cdot \nabla) \psi^i - \br (\psi\cdot \nabla) \d_t \psi^i - \br (\tu\cdot \nabla) \d_t \psi^i - \bq \d_{x_i} \d_t \phi.
    \end{equation*}
    Thus applying Lemma \ref{lemma:elliptic} on $\zeta = \d_t\psi$ with $k=0$, and using inequality \eqref{Linf}, and Lemmas \ref{lemma:timeD}, \ref{lemma:Htwo}, we obtain that
    \begin{align}
        \int_{0}^T\!\! \|D^2 \d_t\psi\|^2\, \dif t  \le C \int_{0}^T\!\! \big\{ \|\mathcal{R}\|^2 + \|\d_t\psi\|^2  \big\}\, \dif t \le M_{\delta}(T). \label{temp:Hthree1}
    \end{align}
    Take $(k,l)=(2,0)$ in Lemma \ref{lemma:damping}. Then by Lemmas \ref{lemma:rE}, \ref{lemma:dpsi}, Corollary \ref{corol:dkdtphi}, and \eqref{temp:Hthree1},
    \begin{align}
        &\sup\limits_{0\le t\le T} \| \sqrt{\chi} D \d^2 \phi\|^2 + \int_{0}^T \!\!\! \big\{ \| \sqrt{\chi} \d_r \d^2 \phi\|^2 + \| \sqrt{\chi} D^3 D_t \phi\|^2  \big\} \dif t 
        \le M_{\delta}(T).\label{temp:Hthree2}
    \end{align}
    Take $(k,l)=(2,0)$ in Lemma \ref{lemma:stokes}. Then by Lemmas \ref{lemma:dpsi}, \ref{lemma:Htwo} and \eqref{temp:Hthree1}--\eqref{temp:Hthree2},
    \begin{align}
        &\int_{0}^T\!\! \big\{ \| \chi D^2 \d^2 \psi\|^2 + \| \chi D \d^2 \phi\|^2 \big\} \dif t 
        \le M_{\delta}(T).\label{temp:Hthree3}
    \end{align}
    Take $(k,l)=(1,1)$ in Lemma \ref{lemma:damping}. Then by Lemmas \ref{lemma:dpsi}, \ref{lemma:Htwo}, Corollary \ref{corol:dkdtphi}, inequalities \eqref{temp:Hthree1}, and \eqref{temp:Hthree3}, we have
    \begin{align}
        &\sup\limits_{0\le t \le T} \| \sqrt{\chi} D^2 \d \phi\|^2 + \int_{0}^T \!\! \Big\{ \big\| \sqrt{\chi} D^2 \d \phi \big\|^2 + \big\| \sqrt{\chi} D^2 \d D_t \phi \big\|^2 \Big\} \dif t 
        \le M_{\delta}(T). \label{temp:Hthree4}
    \end{align}
    Take $(k,l)=(1,1)$ in Lemma \ref{lemma:stokes}. Then by Lemma \ref{lemma:Htwo}, inequalities \eqref{temp:Hthree1} and \eqref{temp:Hthree4},
    \begin{align}
        &\int_{0}^T\!\! \Big\{ \big\| \chi D^3 \d \psi \big\|^2 + \big\| \chi D^2 \d \phi \big\|^2 \Big\}\dif t 
        \le M_{\delta}(T). \label{temp:Hthree5}
    \end{align}
    Take $(k,l)=(0,2)$ in Lemma \ref{lemma:damping}. Then by Lemmas \ref{lemma:dpsi}, \ref{lemma:Htwo}, Corollary \ref{corol:dkdtphi}, inequalities \eqref{temp:Hthree1} and \eqref{temp:Hthree5}, we obtain
    \begin{align}
        &\sup\limits_{0\le t\le T} \big\| D^3 \phi\|^2 + \int_{0}^T\!\! \big\{ \| D^3 \phi\|^2 + \| D^{3} D_t\phi \|^2 \big\} \dif t 
        \le M_{\delta}(T).\label{temp:Hthree6}
    \end{align}
    Take $(k,l)=(0,2)$ in Lemma \ref{lemma:stokes}. Then by Lemma \ref{lemma:Htwo}, inequalities \eqref{temp:Hthree1} and \eqref{temp:Hthree6},
    \begin{align}
        &\int_{0}^T\!\! \big\{ \|D^4 \psi\|^2 + \|D^3 \phi\|^2 \big\} \dif t 
        \le M_{\delta}(T).
    \end{align}
    Finally, from the momentum equation \eqref{Ei}, we see that $\zeta = \psi$ solves system \eqref{elliptic} with $\nu_1 = \mu$, $\nu_2 = \mu+\lambda$, and $\mathcal{R} = \br \{f-\d_t\psi-(u\cdot\nabla)\psi - \bq\nabla \phi\}$. Applying Lemma \ref{lemma:elliptic} on $\zeta= \psi$ with $k=3$, we get $\sup_{0\le t\le T}\|D^3 \psi(t,\cdot)\|^2 \le M_{\delta}(T)$. 
\end{proof}

\begin{proof}[Proof of Theorem \ref{thm:apriori}]
    By \eqref{norm} and Lemmas \ref{lemma:Htwo}--\ref{lemma:Hthree}, there exists a positive constant $C=C(\br,\gamma,K,\mu,\lambda)$ such that for $\delta\in(0,1)$ and $t\in[0,T]$,
    \begin{align}
        N^2(t) \le \frac{C}{\delta^7} \Big\{ N^2(0) + |u_b| N^2(t) + N^3(t) \Big\} + C \delta N^2(t).\label{temp:apriori1}
    \end{align}
    Let $\delta = \frac{1}{4C}$, and set $ \ep = 4^{-9}C^{-8}$. \eqref{continuity} and \eqref{temp:apriori1} implies that if $|u_b| + N(0) \le \ep$, then
    \begin{equation*}
        N^2(t) \le 4^7 C^8 N^2(0) + \frac{3}{4} N^2(t)  \qquad \text{ for all } \ t\in[0,T].
    \end{equation*}
    Therefore $N^2(T) \le 4^8C^8 N^2(0)$, and we conclude the proof of Theorem \ref{thm:apriori}.
\end{proof}

\section*{Declaration of competing interests}
The authors declare that there is no conflict of interests in the present project. The authors also declare that this manuscript has not been previously published, and will not be submitted elsewhere before the decision.

\section*{Funding}
The research of Yucong Huang was supported in part by the UK Engineering and Physical Sciences Research Council Award EP/L015811/1. The research of Shinya Nishibata was supported in part by Grant-in-Aid for Scientific Research (C)(2) 14540200 of the Ministry of Education, Culture, Sports, Science and Technology Science and Technology.

\section*{Declaration of generative AI in scientific writing}
The authors declare that no generative AI is used in the writing of this paper.

\appendix

\section{Appendix}
\setcounter{equation}{0}
\subsection{Hardy's inequality}
\begin{proposition}\label{prop:hardy}
Let $\Omega\vcentcolon=\{ x\in\R^n \,\vert\, |x| \ge 1 \}$ with $n\ge 3$. If $u \in H^1(\Omega)$, then $\frac{u}{|x|}\in L^2(\Omega)$, and the following estimate holds:
\begin{equation*}
\int_{\Omega}\dfrac{|u|^2}{|x|^2}\,\dif x + \int_{\{|x|=1\}}\!\! |u|^2\, \dif S_x \le 2n \int_{\Omega}|\nabla u|^2 \, \dif x.
\end{equation*}
\end{proposition}
\begin{proof}
Fix $R>1$, and let $u\in (\mathcal{C}^{\infty}\cap H^1)(\Omega)$. Since $\nabla(\frac{1}{|x|})=-\frac{x}{|x|^3}$, we have
\begin{align*}
&\int_{\Omega_R} \dfrac{|u|^2}{|x|^2}\, \dif x = - \int_{\Omega_R} |u|^2 \dfrac{x}{|x|} \cdot \nabla\big( \dfrac{1}{|x|} \big)\, \dif x \\
=& \int_{\Omega_R}\!\!\Big\{ 2 u \cdot \nabla u \cdot \dfrac{x}{|x|^2} + (n-1)\dfrac{|u|^2}{|x|^2} \Big\} \, \dif x - \dfrac{1}{R} \int_{\{ |x|=R \}}\!\! |u|^2\, \dif S_x + \int_{\{ |x|=1 \}}\!\! |u|^2\, \dif S_x.
\end{align*}
Thus reorganising the equation we obtain
\begin{align*}
(n-2) \int_{\Omega_R} \dfrac{|u|^2}{|x|^2}\,\dif x + \int_{\{ |x|=1 \}}\!\! |u|^2\, \dif S_x = -2 \int_{\Omega_R} u \cdot \nabla u \cdot \dfrac{x}{|x|^2}\, \dif x + \dfrac{1}{R}\int_{\{ |x|=R \}} \!\! |u|^2\, \dif S_x.
\end{align*}
By Cauchy-Schwartz inequality, it follows that
\begin{align}\label{Atemp1}
&(n-2)\!\!\int_{\Omega_R} \dfrac{|u|^2}{|x|^2}\, \dif x + \int_{\{ |x|=1 \}}\!\! |u|^2\, \dif S_x\nonumber\\
\le& \dfrac{n-2}{2} \int_{\Omega_R} \!\! \dfrac{|u|^2}{|x|^2}\,\dif x + \dfrac{2}{n-2} \int_{\Omega_R}\!\! |\nabla u|^2\, \dif x + \int_{\{ |x|=R \}} \!\! \dfrac{|u|^2}{R}\, \dif S_x
\end{align}
In addition, by divergence theorem and Cauchy-Schwartz inequality, we also have
\begin{align*}
\int_{\{|x|=R\}}\!\! \dfrac{|u|^2}{R}\, \dif S_x - \int_{\{ |x|=1 \}}\!\!\dfrac{|u|^2}{R^2}\,\dif S_x =& \dfrac{1}{R^2}\int_{\Omega_R}\!\! \dd( x |u|^2 ) \, \dif x =\dfrac{1}{R^2}\int_{\Omega_R}\!\! \big\{ n |u|^2 + 2 x u \nabla u  \big\}\, \dif x\\
\le& \dfrac{n+1}{R^2}\int_{\Omega_R}\!\! |u|^2\, \dif x + \int_{\Omega_R}\!\! |\nabla u|^2\, \dif x.  
\end{align*}
Substituting the above inequality
in
\eqref{Atemp1}, it follows that
\begin{align}\label{Atemp2}
\int_{\Omega_R} \dfrac{|u|^2}{|x|^2}\,\dif x + \dfrac{R^2-1}{R^2}\int_{\{|x|=1\}}\!\! |u|^2\, \dif S_x \le 2n \int_{\Omega_R}\!\!|\nabla u|^2 \, \dif x + \dfrac{2(n+1)}{n-2} \int_{\Omega_R}\!\! \dfrac{|u|^2}{R^2} \, \dif x.
\end{align}
By density argument, (\ref{Atemp2}) holds for any $u\in H^1(\Omega)$. We remark that the trace of $u$ on the boundary domain $\{|x|=1\}$ is well defined by the trace embedding theorem. Taking $R\to\infty$ on the resultant inequality, we obtain that
\begin{align*}
\int_{\Omega}\dfrac{|u|^2}{|x|^2}\,\dif x + \int_{\{|x|=1\}}\!\! |u|^2\, \dif S_x \le 2n \int_{\Omega}\!\!|\nabla u|^2 \, \dif x.
\end{align*}
This completes the proof.
\end{proof}

\subsection{Spherical differential operators}\label{append:sph}
In this appendix, we summarize the properties of differential operators in spherical coordinate. First, we give the detailed construction of 
the cut-off functions claimed in \eqref{chi}. Define the function $\tilde{\xi}(\theta)\in \mC^1[0,\pi]$ by
\begin{equation*}
    \tilde{\xi}(\theta) \vcentcolon= \left\{\begin{aligned}
        &0 && \text{if } \ \theta \in [0, \tfrac{\pi}{8})\\ 
        &\tfrac{32}{\pi^2}(\theta-\tfrac{\pi}{8})^2 && \text{if } \ \theta \in [ \tfrac{\pi}{8}, \tfrac{\pi}{4} ),\\
        & 1-\tfrac{32}{\pi^2}(\theta-\tfrac{3\pi}{8})^2 && \text{if } \ \theta \in [ \tfrac{\pi}{4}, \tfrac{3\pi}{8} ),\\
        & 1 && \text{if } \ \theta \in [\tfrac{3\pi}{8}, \tfrac{5\pi}{8}),
    \end{aligned}\right. \quad \left\{\begin{aligned}
        & 1 && \text{if } \ \theta \in [\tfrac{3\pi}{8}, \tfrac{5\pi}{8}),\\
        & 1 - \tfrac{32}{\pi^2}(\theta - \tfrac{5\pi}{8})^2 && \text{if } \ \theta \in [\tfrac{5\pi}{8}, \tfrac{3\pi}{4}),\\
        &\tfrac{32}{\pi^2}(\theta - \tfrac{7\pi}{8})^2 && \text{if } \ \theta \in [\tfrac{3\pi}{4}, \tfrac{7\pi}{8}),\\
        & 0 && \text{if } \ \theta \in [\tfrac{7\pi}{8},\pi].
    \end{aligned}\right.
\end{equation*}
We verify that for $\theta\in [0,\pi]$,
\begin{equation*}
    0\le \tilde{\xi}(\theta) \le 1, \qquad |\tilde{\xi}^{\prime}(\theta)|^2 \le \frac{128}{\pi^2} \xi(\theta), \qquad \supp\big(\tilde{\xi}\big) = \big[\frac{\pi}{8},\frac{7\pi}{8}\big].
\end{equation*}
Using the above properties, we choose a suitable mollifier $\eta\in \mC^{\infty}_{c}[-1,1]$ such that the convolution $\xi\vcentcolon= \tilde{\xi} \ast \eta $ satisfies,
\begin{equation}\label{xi}
    0\le \xi \le 1, \qquad |\xi^{\prime}|^2 \le C \xi, \qquad \supp(\xi) \subseteq \big[ \frac{\pi}{9}, \frac{8\pi}{9} \big], \qquad \sup\limits_{\theta\in[0,\pi]} \big| \xi^{(k)}(\theta) \big| \le C, 
\end{equation}
for $k=1,2,3,4$, where $C>0$ is a generic constant.
Using this, we define the partition of unity functions as
\begin{align*}
    \chi_V(x) \vcentcolon= \xi\Big( \arccos\big( \frac{x_3}{|x|} \big) \Big) \quad \text{ and } \quad \chi_H(x) \vcentcolon= \xi\Big( \arccos\big( \frac{x_2}{|x|} \big) \Big).
\end{align*}
Since $\supp(\xi)=[\frac{\pi}{9},\frac{8\pi}{9}]$, we see that $\chi_V = 0$ for $\arccos(\frac{x_3}{|x|}) \in [0,\frac{\pi}{9}]\cup [\frac{8\pi}{9},\pi]$. The gradient of $\chi_V$ is given by
\begin{align*}
    \nabla \chi_V(x) = \xi^{\prime}\Big( \arccos\big( \frac{x_3}{|x|} \big) \Big) \frac{\hat{e}_1 x_1x_3 + \hat{e}_2 x_2 x_3 - \hat{e}_3(x_1^2 + x_2^2)}{|x| \sqrt{x_1^2+x_2^2}} = \hth_V\, \xi^{\prime}(\theta_V). 
\end{align*}
Thus using \eqref{xi}, we also have for $x\in\Omega$,
\begin{equation*}
    |\nabla \chi_V(x)|^2 \le C \chi_V(x), \qquad \hr\cdot \nabla \chi_V = \hvp_V \cdot \nabla \chi_V = 0.
\end{equation*}
The assertions for $\chi_H(x)$ is shown following the exact same argument above. This proves the existence of functions described in \eqref{chi}. 

In what follows, we use the abbreviated notations $(\theta,\vp) \vcentcolon= (\theta_V,\vp_V)$ or $(\theta_H,\vp_H)$ and $\chi \vcentcolon= \chi_{V}$ or $\chi_H$. Let $(r,\theta,\varphi)\in (0,\infty)\times[0,\pi]\times[0,2\pi)$ be the spherical coordinate. Then the following differential relations for unit basis vectors holds:
\begin{subequations}\label{sph-Dhat}
    \begin{alignat}{4}
    & \d_r \hr = 0,\qquad && \d_\theta \hr = \hth, \qquad && \d_{\varphi} \hr = \hvp \sin\theta,\\
    & \d_r \hth = 0, && \d_\theta \hth = -\hr, && \d_{\varphi} \hth = \hvp \cos\theta,\\
    & \d_r \hvp = 0, && \d_\theta \hvp = 0, && \d_{\varphi} \hvp = -\hr\sin\theta -\hth \cos\theta.
\end{alignat}
\end{subequations}
Let $F\vcentcolon\R^3\to \R$ be a scalar-valued function and $V\vcentcolon\R^3\to \R^3$ be a vector-valued function. Then the gradient, divergence, and Laplacian under spherical coordinate are given by
\begin{subequations}\label{sph-D}
    \begin{gather}
    \nabla F = \hr \d_r F + \hth \frac{1}{r}\d_\theta F + \hvp \frac{1}{r\sin\theta} \d_{\varphi} F,\label{sph-grad}\\
    \dd V = \frac{1}{r^2} \d_r\big(r^2 (\hr\cdot V) \big) + \frac{1}{r\sin\theta}\d_{\theta}\big( (\hth\cdot V) \sin\theta \big) +\frac{1}{r\sin\theta} \d_{\varphi}\big( \hvp\cdot V \big),\label{sph-div}\\
    \Delta F = \frac{1}{r^2} \d_r \big( r^2 \d_r F \big) + \frac{1}{r^2 \sin\theta} \d_\theta\big(\d_\theta F \sin\theta\big) + \frac{1}{r^2\sin^2\theta} \d_{\varphi}^2 F.\label{sph-lapl}
\end{gather}
\end{subequations}

\begin{proposition}\label{prop:dcom}
    Let $\d\in\{\d_\theta,\d_{\vp}\}$ and $D\in \{\d_r, \d_\theta,\d_{\vp} \}$. For an arbitrary non-negative integer $N$, there exists a constant $C=C(N)>0$ such that for an arbitrary scalar functions $F(x)\vcentcolon \R^3 \to \R$, we have
    \begin{gather*}
        \chi\big|[\d^N,\nabla] F\big| \le C \chi \sum_{m=0}^{N-1}|\nabla \d^{m} F|, \qquad \chi \big|[\d^N,\Delta]F\big| \le C\chi 
 \frac{1}{r^2}\sum_{m=1}^{N+1}|\d^{m} F|.
    \end{gather*}
    Moreover, for functions $F(x)\vcentcolon\R^3\to\R$ and $V(x)\vcentcolon \R^3\to\R^3$,
    \begin{align*}
    \chi \big| [D^N,(V\cdot \nabla)]F \big| \le& C \chi \sum_{k=1}^N |D^k F| \sum_{m=0}^{N-k+1} |D^m V|, \\
    \chi \big| [\d^N,\nabla \dd] V \big| \le& C \chi \frac{|V|}{|x|^2} + C \chi \sum_{m=0}^{N-1} |D^2 \d^m V|.
    \end{align*}
\end{proposition}
\begin{proof}
    By product rule, the binomial theorem and \eqref{sph-grad}, it holds that
    \begin{equation*}
        [\d^N, \nabla] F = \sum_{m=0}^{N-1} \binom{N}{m} \Big\{ \big(\d^{N-m}\hr \big) \d_r \d^m F +  \frac{\d^{N-m}\hth}{r}  \d_\theta \d^m F + \d^{N-m} \Big( \frac{\hvp}{r\sin\theta} \Big) \d_{\vp} \d^m F  \Big\}.
    \end{equation*}
    Thus using the differential relations \eqref{sph-Dhat}, we obtain the first inequality. By a similar calculation, we also obtain from \eqref{sph-lapl} that
    \begin{align*}
        [\d^N,\Delta] F = \frac{1}{r^2}\sum_{m=0}^{N-1} \binom{N}{m} \Big\{\d_{\theta} \d^m F \d^{N-m}\Big(\frac{\cos\theta}{\sin\theta}\Big) + \d_{\vp}^2 \d^m F \d^{N-m}\Big(\frac{1}{\sin^2\theta}\Big)\Big\}.  
    \end{align*}
    This yields the second inequality of the proposition. Let $D\in\{\d_r,\d_\theta,\d_\varphi\}$ and $V\vcentcolon\R^3\to\R^3$ be a vector-valued function. We have by the binomial theorem and \eqref{sph-D} that
    \begin{align*}
        &\big[D^N, (V\cdot \nabla )\big] F\\ 
        =& \sum_{k=0}^{N-1} \binom{N}{k} D^{N-k} V \cdot \nabla D^k F + \sum_{m=0}^{N-1} \sum_{k=m+1}^N \binom{N}{k} \binom{k}{m}   D^{N-k} V \cdot D^{k-m} \hr D^m \d_r F \\
        &+ \sum_{m=0}^{N-1} \sum_{k=m+1}^N \binom{N}{k} \binom{k}{m} D^{N-k} V \cdot \Big\{ D^{k-m}\Big(\frac{\hth}{r}\Big) D^m \d_{\theta} F + D^{k-m} \Big(\frac{\hvp}{r\sin\theta}\Big) D^m \d_\vp F \Big\}. 
    \end{align*}
    Multiplying both sides by $\chi$ and taking the absolute value, we get
    \begin{align*}
        \big|[D^N,(V\cdot \nabla)]F\big| \le & C\chi \Big\{ \sum_{k=0}^{N-1} |D^{N-k}V| |D^{k+1}F|  + \sum_{k=0}^{N-1} |D^{k+1}F| \sum_{m=k+1}^N |D^{N-m}V| \Big\}\\
        =& C \chi \sum_{k=1}^N |D^k F| \sum_{m=0}^{N-k+1} |D^m V|.
    \end{align*}
    By \eqref{sph-grad}--\eqref{sph-div}, we have 
    \begin{align*}
        &\nabla \dd V \\
        =& \hr (\hr\cdot  \d_r^2 V) + \hth \frac{\d_r\d_\theta (\hr\cdot V)}{r} + \hvp\frac{\d_r\d_\vp (\hr\cdot V)}{r\sin\theta}  + 2\hr \Big( \frac{\d_r(\hr\cdot V)}{r} - \frac{\hr\cdot V}{r^2} \Big) + \hth \frac{2\d_\theta(\hr\cdot V)}{r^2} \\
        & + \hvp \frac{2 \d_{\vp}(\hr\cdot V)}{r^2 \sin\theta} + \hr \Big( \frac{\d_r\d_\theta(\hth\cdot V)}{r} - \frac{\d_\theta (\hth\cdot V)}{r^2} \Big) + \hth \frac{\d_\theta^2 (\hth \cdot V)}{r^2} + \hvp \frac{\d_\vp\d_\theta (\theta \cdot V)}{r^2 \sin\theta} \\
        & + \hr \Big( \frac{\d_r(\hth\cdot V) \cos\theta}{r\sin\theta} - \frac{(\hth\cdot V) \cos\theta}{r^2 \sin\theta} \Big) + \hth \d_\theta \Big( \frac{(\hth\cdot V)\cos\theta}{r^2 \sin\theta} \Big) + \hvp \frac{\d_\vp(\hth\cdot V) \cos\theta}{r^2 \sin^2\theta} \\
        & + \hr \Big( \frac{\d_r \d_{\vp} (\hvp\cdot V)}{r\sin\theta} - \frac{\d_\vp (\hvp\cdot V)}{r^2 \sin\theta} \Big) + \hth \d_\theta\Big( \frac{\d_\vp(\hvp\cdot V)}{r^2 \sin\theta} \Big) + \hvp \frac{\d_{\vp}^2 (\hvp \cdot V)}{r^2 \sin^2\theta}.
    \end{align*}
    Computing the commutator $[\d^k, \nabla \dd]$ using the above expression, we obtain the last inequality in the proposition.
\end{proof}

\begin{proposition}\label{prop:rrcancel}
    Let $\d\in\{\d_\theta,\d_{\vp}\}$ and $D\in \{\d_r, \d_\theta,\d_{\vp} \}$. Then for $k,l \in \mathbb{N}\cup\{0\}$,
    there exists a constant $C=C(k,l)>0$ such that for an arbitrary vector-valued function $V(x)\vcentcolon \R^3 \to \R^3$, it holds that
\begin{align*}
    \int_{\Omega} \chi \big| \d_r^{l} \d^k (\hr\cdot \Delta V - \d_r \dd V) \big|^2 \dif x 
    \le& C \int_{\Omega} \chi \sum_{i=1}^{l+1} \sum_{j=0}^{k+1} \big|  D^{i} \d^{j} V \big|^2 \dif x .  
\end{align*}
\end{proposition}
\begin{proof}
Using differential relations \eqref{sph-Dhat} and \eqref{sph-D}, we get
\begin{align*}
    \hr\cdot\Delta V 
    =& \hr\cdot \d_r^2 V + \frac{2}{r} \hr\cdot \d_r V + \frac{\cot\theta}{r^2} \hr\cdot \d_\theta V + \frac{1}{r^2} \hr \cdot \d_\theta^2 V + \frac{1}{r^2 \sin^2\theta} \hr \cdot \d_{\vp}^2 V,\\
    \d_r \dd V =& \hr \cdot \d_r^2 V + \frac{\hth\cdot \d_{\theta} \d_r V }{r} + \frac{\hvp \cdot \d_{\vp} \d_r V}{r\sin\theta}  - \frac{\hth\cdot \d_{\theta} V}{r^2}  - \frac{\hr\cdot V}{r^2}  - \frac{\cot\theta}{r^2} \hth\cdot V. 
\end{align*}
Taking the difference between the above two equations, we see that the term $\hr\cdot \d_r^2 V$ cancels. Using this and Proposition \ref{prop:hardy}, we obtain the desired estimate.
\end{proof}

\end{document}